\documentclass[a4paper,12pt]{article}
\usepackage[cp1251]{inputenc}            
\usepackage[english]{babel}
\usepackage[left=2.5cm,right=1.5cm,top=2cm,bottom=2cm]{geometry} 
\usepackage{amsmath,amssymb,amsthm}
\usepackage{envmath}
\usepackage{mathrsfs}
\usepackage{caption2}
\usepackage{floatflt}
\usepackage[multiple]{footmisc}

\newtheorem{theorem}{Theorem}
\newtheorem{lemma}{Lemma}
\newtheorem{example}{Example}

\newtheorem{definition}{Definition}

\newtheorem{corollary}{Corollary}

\newtheorem{remark}{Remark}

\addto\captionsrussian{}

\allowdisplaybreaks[1]

\begin{document}







\title{The limit joint distributions of some statistics used in testing the quality of random number generators}

\author{Savelov M.P.\footnote{Lomonosov Moscow State University, maksim.savelov@math.msu.ru}}

\maketitle

\begin{abstract}
The limit joint distribution of statistics that are generalizations of some statistics from the NIST STS, TestU01, and other packages is found under the following hypotheses $H_0$ and $H_1$. Hypothesis $H_0$ states that the tested sequence is a sequence of independent random vectors with a known distribution, and the simple alternative hypothesis $H_1$ converges in some sense to $H_0$ with increasing sample size. In addition, an analogue of the Berry-Esseen inequality is obtained for the statistics under consideration, and conditions for their asymptotic independence are found.
\end{abstract}
\textbf{Keywords:} joint distribution of statistics, NIST STS, TestU01, summing statistics, long-block statistics, short-block statistics, quadratic statistics, goodness-of-fit tests, asymptotic independence, Berry-Esseen-type inequality

\section{Introduction}

Testing the quality of random number generators is an important practical problem. Of particular interest is testing the quality of random binary sequence generators, i.e. testing the hypothesis $H_0^{Bern(\frac12)}$ that the tested sequence $\varepsilon_1, \ldots, \varepsilon_n$ is a sequence of independent Bernoulli trials with parameter $\frac12$. Another important problem is testing the quality of random number generators from the segment $[0,1]$, which tests the hypothesis $H_0^{U[0,1]}$ that the tested sequence $\varepsilon_1, \ldots, \varepsilon_n$ is a sequence of independent random variables with uniform distribution on the segment $[0,1]$. There are a number of tools for testing $H_0^{Bern(\frac12)}$ and $H_0^{U[0,1]}$: NIST STS \cite{Nist}, TestU01 \cite{LecyerTestU01math, TestU01ShortGuide} and Crypt-X packages, Diehard tests, D. Knuth tests \cite{Knuth}, etc. Each of them contains a set of several statistical tests. When applying them, the question of the limit joint distribution of their statistics naturally arises. There are a number of paper  devoted to the limit joint distributions of statistics used in testing $H_0^{Bern(\frac12)}$, see, e.g., \cite{Zubkov, Voloshko3, KharinZubkov, Savelov11, Savelov12, Savelov13} (for earlier publications, see the bibliography in \cite{Savelov12}). However, little attention was paid to the limit joint distributions of statistics used in testing the quality of random number generators with a uniform distribution on the segment $[0,1]$. The results of this paper are applicable both to statistics used in testing $H_0^{Bern(\frac12)}$ and to statistics used in testing $H_0^{U[0,1]}$.

The most popular packages for testing the hypothesis $H_0^{Bern(\frac12)}$ are NIST STS and TestU01. In \cite{Savelov11} and \cite{Savelov12} in the case, when the elements of the tested sequence take values in a fixed finite set $\{0,1,\ldots,R-1\}$, four types of statistics were defined: summing, long-block, short-block, and quadratic statistics. These statistics are generalizations of several statistics of tests from NIST STS, TestU01, and other packages. The limit joint distribution of the four types of statistics is found in \cite{Savelov12}. In the present paper, summing, long-block, short-block, and quadratic statistics are defined in the general case (when the distribution of sample elements is not necessarily discrete). The limit joint distribution of these statistics is found in the case where either the hypothesis $H_0$ (that states that the elements of the tested sequence are independent and have a fixed pre-assigned distribution) is true, or a simple alternative $H_1$ is true, where $H_1$ is in some sense "close" to $H_0$ (see Theorem \ref{thOsovmRaspr3stat}, Corollary \ref{thOsovmRaspr3stat}, and Remark \ref{zamechProH1BlizkH0}). This result generalizes some of the results of \cite{Savelov11} and \cite{Savelov12} (since it concerns not only discrete samples), and also allows us to generalize the most important results of \cite{Savelov12} and \cite{Savelov13} (based on \cite{Savelov11}), i.e. obtain an analogue of the Berry-Esseen inequality for summing, long-block, and short-block statistics (see Theorem \ref{th2Berry}) and necessary and sufficient conditions for asymptotic independence of statistics of four types (see Theorem \ref{thMainAsNez}). In addition, in this paper we obtain conditions under which the summing/long-block/short-block statistics tend to infinity (see Lemma \ref{LemmaProNeBlizkH1} and Remark \ref{ProLemmyProNeblizkH1}).

The most complete set of tests for testing the hypothesis $H_0^{U[0,1]}$ is contained in the TestU01 package (see the BigCrush battery). The statistics of some of them are statistics of one of the four types considered --- see examples in Section \ref{sectionExamples}. The limit joint distribution of such statistics cannot be found using the results of \cite{Savelov11, Savelov12, Savelov13}, which deal with the case of a discrete sample, but can be found using Theorem \ref{thOsovmRaspr3stat} and Corollary \ref{thOsovmRaspr3stat}.

\section{Basic Definitions}

The definitions in this section are almost similar to the corresponding definitions in \cite{Savelov11} and \cite{Savelov12}.

Let $d$ be a non-negative integer and let $\varepsilon_1, \ldots, \varepsilon_n$ be a sequence of random variables, each taking values in the set $\mathfrak{X}$. For simplicity, we can assume that $\mathfrak{X} \subset \mathbb{R}^d$, although a more general case can be considered. Of greatest interest are the cases where $\mathfrak{X} = [0,1]$ or $\mathfrak{X} = \{0,1\}$. By $H_0$ we denote the simple hypothesis that the sequence $\varepsilon_1, \ldots, \varepsilon_n$ consists of independent random variables with a fixed distribution $\mathbf{P}_0$, i.e. $\mathbf{P}_{H_0}(\varepsilon_i \in B) = \mathbf{P}_0(B)$ for all measurable sets $B$ from the corresponding sigma-algebra (in the case of $\mathfrak{X} \subset \mathbb{R}^d$, the Borel sigma-algebra is considered).

Let ${N_{lb}}$ and $L_{sb}$ be natural numbers (they are the parameters of the statistics $T_{lb}$ and $T_{sb}$ constructed below). We will be interested in the case when ${N_{lb}}$ and $L_{sb}$ are fixed and $n$ tends to infinity. Put $${L_{lb}} = \left[\frac{n}{N_{lb}}\right], \ \ \ {N_{sb}} = \left[\frac{n}{L_{sb}}\right].$$ If we discard the last $n - {N_{lb}}{L_{lb}}$ elements from the original sequence of random variables $\varepsilon_1, \ldots, \varepsilon_n$, then the remaining elements are splitted into ${N_{lb}}$ ''long'' disjoint blocks of length ${L_{lb}}$: the first block consists of random variables $\varepsilon_1, \ldots, \varepsilon_{L_{lb}}$, the second  --- of random variables $\varepsilon_{{L_{lb}}+1}, \ldots, \varepsilon_{2{L_{lb}}}$, etc. Similarly, we can discard the last $n - L_{sb}{N_{sb}}$ elements and split the remaining ones into ${N_{sb}}$ ''short'' blocks of length $L_{sb}$.

The subscript "sb" denotes that the corresponding value refers to the statistics constructed over short blocks. Similarly, the subscript "lb" is a shorthand for "long blocks."

Further, the Cartesian product $\underbrace{\mathfrak{X} \times \mathfrak{X} \times \ldots \times \mathfrak{X}}_{\text{$k$ times}}$ will be denoted by $\mathfrak{X}^{k}$. We fix natural numbers ${m_{sum}}, m_{lb}, K_{sb}$ and three functions
$f_{sum}: \mathfrak{X}^{m_{sum}} \to \mathbb{R}$,
$f_{lb}: \mathfrak{X}^{m_{lb}} \to \mathbb{R}$ and
$f_{sb}: \mathfrak{X}^{L_{sb}} \to \mathbb{R}$. Partition the set of values of $f_{sb}$ into ${K_{sb}}+1$ nonempty disjoint subsets: $f_{sb}(\mathfrak{X}^{L_{sb}}) = \bigsqcup_{i=0}^{K_{sb}} \alpha_{sb}(i)$. We will use the following notation. The quantity $\mathbf{E}_{H_0} f_{sum}(\varepsilon_1, \ldots, \varepsilon_{{m_{sum}}})$ is the expectation of $f_{sum}(\varepsilon_1, \ldots, \varepsilon_{{m_{sum}}})$, calculated under the assumption that the hypothesis $H_0$ is true. The variance $\mathbf{D}_{H_0} f_{sum}(\varepsilon_1, \ldots, \varepsilon_{{m_{sum}}})$, etc., are interpreted similarly. We will be interested in the case when the numbers $m_{sum}, m_{lb}, N_{lb}, L_{sb}, K_{sb}$, the functions $f_{sum}, f_{lb}, f_{sb}$ and the sets $\alpha_{sb}(0), \ldots, \alpha_{sb}(K_{sb})$ (as well as the distribution $\mathbf{P}_{H_0}$) are fixed, do not depend on $n$, and $n$ tends to infinity.

Put\begin{gather}
E_{sum} = \mathbf{E}_{H_0} f_{sum}(\varepsilon_1, \ldots, \varepsilon_{{m_{sum}}}), \nonumber \\
\sigma_{sum}^2 = \mathbf{D}_{H_0} f_{sum}(\varepsilon_1, \ldots, \varepsilon_{{m_{sum}}}) + \nonumber \\
+2\sum_{i=2}^{{m_{sum}}} \mathrm{cov}_{H_0}\big( f_{sum}(\varepsilon_{i}, \ldots, \varepsilon_{i+{m_{sum}}-1}), f_{sum}(\varepsilon_{1}, \ldots, \varepsilon_{{m_{sum}}}) \big)\label{dlyaSigmaSquared}
\end{gather}
and $\sigma_{sum} \ge 0$. Note that the right-hand side of the formula \eqref{dlyaSigmaSquared} is non-negative, since it is the variance of the limit distribution in the central limit theorem for $({m_{sum}}-1)$-dependent random variables $f_{sum}(\varepsilon_i, \ldots, \varepsilon_{i+{m_{sum}}-1})$, $i = 1, \ldots, n-m_{sum}+1$, as $n \to \infty$.
\begin{definition}
If $\sigma_{sum} \in (0, \infty)$, then a statistic of the form
\begin{equation}
T_{sum} = \frac{ \sum_{i=1}^{n-m_{sum}+1} \big( f_{sum}(\varepsilon_i, \ldots, \varepsilon_{i+{m_{sum}}-1}) - E_{sum}\big)}{\sigma_{sum} \sqrt{n - m_{sum} +1}}
\label{defTsum}
\end{equation}
is called the {\bfseries summing} statistic.
\end{definition}

Put \begin{gather}
W_k = \sum_{j={L_{lb}}(k-1)+1}^{{L_{lb}}k-{m_{lb}}+1} f_{lb}(\varepsilon_{j}, \ldots, \varepsilon_{j+{m_{lb}}-1}), \quad 1 \le k \le {N_{lb}}, \nonumber \\
E_{lb}= \mathbf{E}_{H_0} f_{lb}(\varepsilon_1, \ldots, \varepsilon_{{m_{lb}}}), \nonumber \\
\sigma_{lb}^2 = \mathbf{D}_{H_0} f_{lb}(\varepsilon_1, \ldots, \varepsilon_{{m_{lb}}}) + 2\sum_{i=2}^{{m_{lb}}} \mathrm{cov}_{H_0}\big( f_{lb}(\varepsilon_{i}, \ldots, \varepsilon_{i+{m_{lb}}-1}), f_{lb}(\varepsilon_{1}, \ldots, \varepsilon_{{m_{lb}}}) \big), \nonumber
\end{gather}
where it is assumed that $\sigma_{lb} \ge 0$.
Note that $\mathbf{E}_{H_0}W_k = ({L_{lb}}-{m_{lb}}+1)E_{lb}$, $ 1 \le k \le N_{lb}$.
\begin{definition}
Let $\sigma_{lb} \in (0, \infty)$. A statistic of the form $$T_{lb} = \frac{\sum_{k=1}^{{N_{lb}}} \big(W_k - ({L_{lb}}-{m_{lb}}+1) E_{lb}\big)^2}{{L_{lb}} \sigma_{lb}^2}$$
is called the {\bfseries long-block} statistic, constructed over ${N_{lb}}$ blocks.
\end{definition}

For $0 \le j \le {K_{sb}}$ put $$w(j)= \sum_{i=1}^{{N_{sb}}} I_{f_{sb}(\varepsilon_{L_{sb}(i-1)+1}, \ldots, \varepsilon_{L_{sb}i}) \in \alpha_{sb}(j) },$$
$$E_{sb}(j) = \mathbf{P}_{H_0}\big({f_{sb}(\varepsilon_{1}, \ldots, \varepsilon_{L_{sb}}) \in \alpha_{sb}(j) }\big).$$
Note that $\mathbf{E}_{H_0}w(j) = N_{sb} E_{sb}(j), \ 0 \le j \le {K_{sb}}$.
\begin{definition}
Let $E_{sb}(j) > 0$ for all $0 \le j \le {K_{sb}}$. A statistic of the form
\begin{equation}
T_{sb} = \sum_{j=0}^{{K_{sb}}} \frac{(w(j) - N_{sb} E_{sb}(j))^2}{ N_{sb} E_{sb}(j)} \nonumber
\end{equation}
is called the {\bfseries short-block} statistic constructed over blocks of length $L_{sb}$.
\end{definition}

Next, for each $1 \le q \le Q_{quad}$ fix a natural number $m_{sum}^{[q]}$ and a function $f_{sum}^{[q]}: \mathfrak{X}^{m_{sum}^{[q]}} \to \mathbb{R}$ and set
\begin{equation}
T_{sum}^{[q]} = \frac{ \sum_{i=1}^{n-m_{sum}^{[q]}+1} \Big(f_{sum}^{[q]}\big(\varepsilon_i, \ldots, \varepsilon_{i+m_{sum}^{[q]}-1}\big) - E_{sum}^{[q]}\Big)}{\sigma_{sum}^{[q]} \sqrt{n - m_{sum}^{[q]}+1}}, \nonumber
\end{equation}
where the quantities $E_{sum}^{[q]}$ and $\sigma_{sum}^{[q]}$ are defined in the same way as the quantities $E_{sum}$ and $\sigma_{sum}$ from \eqref{defTsum}. The statistics $T_{sum}^{[1]}, \ldots, T_{sum}^{[Q_{quad}]}$ are summing statistics.
\begin{definition}
Consider a positive integer $\tau_{quad}$ and real numbers $d_{quad}(i,q)$ ($1 \le i \le \tau_{quad}, 1 \le q \le Q_{quad}$). A statistic of the form
\begin{equation}
T_{quad} = \sum_{i=1}^{\tau_{quad}} \Bigg(\sum_{q=1}^{Q_{quad}} d_{quad}(i,q) T_{sum}^{[q]}\Bigg)^2 \nonumber
\end{equation}
is called the {\bfseries quadratic} statistic constructed from the statistics $T_{sum}^{[1]}, \ldots, T_{sum}^{[Q_{quad}]}$.
\end{definition}
The index ''quad'' corresponds to the word ''quadratic''. It is easy to see (cf. Section 2 of \cite{Savelov12}) that quadratic statistics are those and only those statistics that are nonnegative definite quadratic forms of the corresponding summing statistics.

As noted above, in the case when the elements of the test sequence take values in a fixed finite set $\{0,1,\ldots,R-1\}$, summing, long-block, short-block, and quadratic statistics were defined in \cite{Savelov11} and \cite{Savelov12}. Furthermore, in \cite{Voloshko3}, under the assumption that the test sequence is binary, statistics similar to summing statistics were considered.

We will call the numbers $m_{sum}, m_{lb}, N_{lb}, L_{sb}, K_{sb}$, the functions $f_{sum}, f_{lb}, f_{sb}$ and the sets $\alpha_{sb}(0), \ldots, \alpha_{sb}(K_{sb})$ the parameters of the statistics $T_{sum}, T_{lb}, T_{sb}$. If these parameters are fixed (independent of $n$) and $H_0$ is true, then by virtue of the considerations from Section 2 \cite{Savelov11}
\begin{equation}
T_{sum} \overset{d}{\to} \mathcal{N}(0,1),\quad T_{lb} \overset{d}{\to} \mathcal{\chi}_{N_{lb}}^2, \quad T_{sb} \overset{d}{\to} \mathcal{\chi}_{K_{sb}}^2, \quad n \to \infty.\label{NotJointDistrOfEvery}
\end{equation}

The numbers $Q_{quad}, \tau_{quad}$, and $d_{quad}(i,q)$ ($1 \le i \le \tau_{quad}, 1 \le q \le Q_{quad}$), as well as the parameters of the statistics $T_{sum}^{[1]}, \ldots, T_{sum}^{[Q_{quad}]}$, are called the parameters of the statistics $T_{quad}$. Further we will assume that they don't depend on $n$ and are fixed.

\section{Examples of summing, long-block, short-block, and quadratic statistics} \label{sectionExamples}

If a sequence of random vectors $\varkappa_n$ converges to zero in probability, then we write  $\varkappa_n = o_{P}(1)$, $n \to \infty$. If the difference between two statistics is $o_{P}(1)$, then we say that they coincide up to $o_{P}(1)$.

We will show how summing, long-block, short-block, and quadratic statistics are related to the statistics of well-known tests used for discrete samples. Note that the statistic of the classical chi-square test is both short-block and quadratic statistic. Moreover, in the case where the elements of the tested sequence take values in the set $\{0,\ldots,R-1\}$ (i.e., if $\mathfrak{X} = \{0,\ldots,R-1\}$), the statistics of "run test" and "poker test" from D. Knuth's battery of tests (see \cite[\S 3.3.2]{Knuth}) are special cases of statistics of the form $T_{sb}$.

Regarding the tests of the NIST STS, it was proved in \cite{Savelov11} that if $\mathfrak{X} = \{0,1\}$ (and the hypothesis $H_0^{Bern(\frac12)}$ is tested), then the statistics $T_{fr}$, $T_{templ}$ of ''Frequency Test within a Block'' and ''Non-overlapping Template Matching Test'' are long-block statistics, the statistics $T_{longrun}, T_{matrix}, T_{lincompl}$ of ''Test for the Longest Run of Ones in a Block'', ''Binary Matrix Rank Test'' and ''Linear Complexity Test'' are short-block statistics (taking into account the caveat from Example 7 \cite{Savelov11}), the statistic $T_{mon}$ of ''Monobit Test'' is a summing statistic, and the statistic $T_{runs}$ of the ''Runs Test'' coincides (under certain additional conditions) with the corresponding summing statistic up to $o_P(1)$ (see Lemma 7 \cite{Savelov11}). Moreover, if $\mathfrak{X} = \{0,1\}$, then by Lemma 6 \cite{Savelov12} the statistics $T_{serial1}, T_{serial2}$ of the ''Serial Test''  and the statistic $T_{entr}$ of the ''Approximate Entropy Test'' are such that each of them, in a wide class of cases, coincides up to $o_P(1)$ with some quadratic statistic. The connections linking the summing, long-block, short-block, and quadratic statistics to those of the Crypt-X and FIPS test suites were discussed in \cite{Savelov11}.

Let us show how the summing, long-block, short-block, and quadratic statistics are related to the statistics of the known tests used in the case where $\mathfrak{X} = [0,1]$ to test the hypothesis $H_0^{U[0,1]}$. Accordingly, in the following examples we will assume that $\mathfrak{X} = [0,1]$.

\begin{example} \label{primerWeightDistrib}
If $0 \le \alpha < \beta \le 1$ and $f_{sb}(\theta_1,\ldots,\theta_{L_{sb}}) = \sum_{i=1}^{L_{sb}}I_{\theta_i \in [\alpha, \beta)}$, then the corresponding short-block statistic $T_{sb}$ coincides with the statistic of the ''svaria\_WeightDistrib'' test of the TestU01 \cite{TestU01ShortGuide} package.
\end{example}

\begin{example} \label{primerSampleCorr}

The ''svaria\_SampleCorr'' test from the TestU01 package uses the statistic
$T_{SampleCorr} = \sqrt{\frac{12}{n-k}} \sum_{i=1}^{n-k} \big( \varepsilon_i \varepsilon_{i+k} - \frac{1}4 \big)$, and it is assumed that $T_{SampleCorr} \overset{d}{\to} \mathcal{N}(0,1)$ as $n \to \infty$. If $m_{sum} \ge 2$ and $f_{sum}(\theta_1,\ldots, \theta_{m_{sum}}) = \theta_1 \theta_{m_{sum}}$, then $$T_{sum} = \frac{12}{\sqrt{13(n-m_{sum}+1)}} \sum_{i=1}^{n-m_{sum}+1} \bigg( \varepsilon_i \varepsilon_{i+m_{sum}-1} - \frac{1}4 \bigg)$$
(cf. \cite[p.245]{SchmidtTaylor}). If $m_{sum} = k+1$, then $T_{sum}$ coincides with $T_{SampleCorr}$ up to normalization, and due to \eqref{NotJointDistrOfEvery}, it is more correct to use $T_{sum}$ instead of $T_{SampleCorr}$ as an approximation to $\mathcal{N}(0,1)$. Note that the description of the ''svaria\_SampleCorr'' test in \cite[p.117]{TestU01ShortGuide} includes the parameter $N$, but the test is only applied when $N=1$ (see \cite[p.146 and p.150]{TestU01ShortGuide}).
\end{example}

\begin{example}
Let us number the elements of the permutation group of the elements of the set $\{1,\ldots, L_{sb}\}$ in an arbitrary but fixed way. Let $\theta_1, \ldots, \theta_{L_{sb}}$ be numbers from the segment $[0,1]$. If they are all pairwise distinct, then the vector $(\theta_1, \ldots, \theta_{L_{sb}})$ naturally corresponds to one of the $L_{sb}!$ elements of the permutation group under consideration. Denote the number of this element by $f_{sb}(\theta_1, \ldots, \theta_{L_{sb}})$. If among the numbers $\theta_1, \ldots, \theta_{L_{sb}}$ there are identical ones, one can extend the function $f_{sb}(\theta_1, \ldots, \theta_{L_{sb}})$ arbitrarily (if the elements of the tested sequence are independent and have a continuous distribution, then the probability of their coincidence is zero). The corresponding short-block statistic $T_{sb}$ coincides with the statistic of the ''sknuth\_Permutation'' test of the TestU01 package and the statistic of "permutation test" from \cite[\S 3.3.2]{Knuth}.
\end{example}

If $\mathfrak{X} = [0,1]$, then, as noted in \cite[\S 3.3.2]{Knuth}, to test the hypothesis $H_0^{U[0,1]}$, one can consider a new sequence $[R\varepsilon_1], [R\varepsilon_2], \ldots$, where $R$ is a ''not too small'' non-negative integer, and apply to it the test used for $\mathfrak{X} = \{0,\ldots,R-1\}$. For example, one can apply the poker test and the run test from \cite[\S 3.3.2]{Knuth}, which use short-block statistics, or the serial test from \cite{Good1,Good2} (see also \cite{Marsaglia2005}), which uses quadratic statistics, etc. Moreover, if $\mathfrak{X} = [0,1]$, then one can apply tests to the sequence $\varepsilon_1, \varepsilon_2, \ldots$ that are based on the assumption that the tested sequence is binary, if one first goes from $\varepsilon_i$ to their binary notation. Here are some examples. To describe them, we need the function $g_{bits}:[0,1]\to\{0,1\}^{r_{bits}}$, which assigns to each number from the segment $[0,1]$ some $r_{bits}$ bits from its binary notation (see, e.g., \cite[\S 3.1]{LEcyer1996}).

If $x \in \mathbb{R}$, denote by $\lceil x \rceil$ the smallest integer that is not less than $x$. In the following examples, we assume that $\mathfrak{X} = [0,1]$.

\begin{example} \label{primerTmatrixRank}
Let $V_1, V_2 \in \mathbb{N}, L_{sb} = V_1 \cdot \big\lceil \frac{V_2}{r_{bits}} \big\rceil$ and $\theta_i \in [0,1]$ for $1 \le i \le L_{sb}$. Construct a matrix $A$ of size $V_1 \times V_2$ as follows. Applying the function $g_{bits}$ to each of the values $\theta_1, \ldots, \theta_{\big\lceil \frac{V_2}{r_{bits}} \big\rceil}$, we obtain $r_{bits} \cdot \big\lceil \frac{V_2}{r_{bits}} \big\rceil$ random bits, the first $V_2$ of which we write in the first row of the matrix $A$. Applying the function $g_{bits}$ to each of the quantities $\theta_{\big\lceil \frac{V_2}{r_{bits}} \big\rceil + 1}, \ldots, \theta_{2\big\lceil \frac{V_2}{r_{bits}} \big\rceil}$, we obtain another $r_{bits} \cdot \big\lceil \frac{V_2}{r_{bits}} \big\rceil$ bits, the first $V_2$ of which we write in the second row of the matrix $A$, and so on. As a result, we construct a matrix $A$ of size $V_1 \times V_2$. Let $f_{sb}(\theta_1, \ldots, \theta_{L_{sb}})$ be the rank of this matrix over the field $GF(2)$. The corresponding short-block statistic $T_{MatrixRank}$ is used in the ''smarsa\_MatrixRank'' test of the TestU01 package.
\end{example}

\begin{example} \label{primerHammingWeight2}
Let $m_{lb} = 1$ and $f_{lb}(\theta) = \sum_{i \in g_{bits}(\theta)}i$ for $\theta \in [0,1]$. The corresponding long-block statistic $T_{lb}$ is analogous to that of the ''sstring\_HammingWeight2'' test of the TestU01 package.
\end{example}

\begin{example} \label{primersmarsaSerialOver}
The ''smarsa\_SerialOver'' test of the TestU01 package is based on the following idea (see \cite{Altman1988}): from a sequence $\varepsilon_1, \ldots, \varepsilon_n$, a sequence of binary vectors $g_{bits}(\varepsilon_1), \ldots, g_{bits}(\varepsilon_n)$ is constructed, each of which takes one of $R = 2^{r_{bits}}$ values (which is equivalent to the case when $\mathfrak{X} = \{0, \ldots, R-1\}$), after which the Serial Test from \cite{Good1,Good2} is applied to the resulting sequence (see also \cite{Marsaglia2005}). Below, in Subsection \ref{oStatistTSo}, the explicit form of the statistic $T_{SO}$ is given and it is proven to be identical to the corresponding quadratic statistic up to $o_P(1)$ when $H_0^{U[0,1]}$ is true (see Lemma \ref{lemmaProTSo}). Note that the NIST package's "Serial Test" uses a statistic similar to the $T_{SO}$ statistic, which is identical to the corresponding quadratic statistic up to $o_P(1)$ when $H_0^{Bern(\frac12)}$ is true (see Lemma 6 \cite{Savelov12}).
\end{example}

Regarding the TestU01 package, we note that it has three batteries of tests for testing $H_0^{U[0,1]}$: SmallCrush, Crush, and BigCrush. The first of them allows to applicate different tests to the same test sequence (see SmallCrushFile in \cite[p.143]{TestU01ShortGuide}), while the BigCrush battery is the biggest one.

Along with Examples \ref{primerWeightDistrib}--\ref{primersmarsaSerialOver}, we can provide other examples of statistics of four types (summing, long-block, short-block, and quadratic) used in various well-known statistical tests.

\section{Main results} \label{MainResults}

\subsection{The limit joint distribution of statistics of 4 types}

First, we formulate a statement that is somewhat less cumbersome and less general than Theorem \ref{th2} below.

Let $N, s, h$ be natural numbers and
\begin{equation}
L = \bigg[\frac{n}{N}\bigg].\nonumber
\end{equation}
We are interested in the case where $N, s$, and $h$ are fixed, and $n$ and $L$ tend to infinity.

Put $$\underline{s} = h + \max({m_{sum}},m_{lb}) -1.$$
If an integer $x$ is a divisor of an integer $y$, then we write $x | y$. Suppose that $N_{lb} | N$ and ${L_{sb}} | h$. Define the function $f^{(\underline{s})}: \mathfrak{X}^{\underline{s}} \to \mathbb{R}^{{K_{sb}}+3}$. Put
\begin{gather}
f_j^{({\underline{s}})}(\theta_{1}, \ldots, \theta_{{\underline{s}}}) = \begin{cases}
\frac{1}{\sqrt{E_{sb}(j)}}\sum_{u=1}^{h / L_{sb}}I_{f_{sb}(\theta_{{L_{sb}}(u-1)+1}, \ldots, \theta_{{L_{sb}}u})\in \alpha_{sb}(j) } ,&\text{if $0 \le j \le {K_{sb}}$,}\\
\sigma_{lb}^{-1} \cdot \sum_{u=1}^{h} f_{lb}(\theta_{u}, \ldots, \theta_{u+{m_{lb}}-1}),&\text{if $j = {K_{sb}}+1$,}\\
\sigma_{sum}^{-1} \cdot \sum_{u=1}^{h} f_{sum}(\theta_u, \ldots, \theta_{u+{m_{sum}}-1}),&\text{if $j = {K_{sb}}+2$;}
\end{cases}\nonumber \\
f^{(\underline{s}) }(\theta_{1}, \ldots, \theta_{{\underline{s}}}) = \Big(f_0^{({\underline{s}})}(\theta_{1}, \ldots, \theta_{{\underline{s}}}), f_1^{({\underline{s}})}(\theta_{1}, \ldots, \theta_{{\underline{s}}}), \ldots, f_{{K_{sb}}+2}^{({\underline{s}})}(\theta_{1}, \ldots, \theta_{{\underline{s}}}) \Big).\nonumber
\end{gather}
Let $s \ge {\underline{s}}$. For each $(\theta_{1}, \ldots, \theta_{s}) \in \mathfrak{X}^s$ put
\begin{gather}
f_j^{(s)}(\theta_{1}, \ldots, \theta_{s}) = f_j^{({\underline{s}})}(\theta_{1}, \ldots, \theta_{{\underline{s}}}), \quad 0 \le j \le {K_{sb}}+2, \nonumber \\
f^{(s)}(\theta_{1}, \ldots, \theta_{s}) = f^{({\underline{s}})}(\theta_{1}, \ldots, \theta_{{\underline{s}}}).\nonumber
\end{gather}

We need the following definitions. If $\vec{\varkappa} = (\varkappa_1, \ldots, \varkappa_{d_1})$ and $\vec{\varkappa}^* = (\varkappa_1^*, \ldots, \varkappa_{d_2}^*)$ are random vectors, then their covariance is defined as the matrix ${\rm cov}(\vec{\varkappa}, \vec{\varkappa}^*) = \mathbf{E} (\vec{\varkappa} - \mathbf{E} \vec{\varkappa})^{\top}(\vec{\varkappa}^* - \mathbf{E} \vec{\varkappa}^*)$ of size $d_1 \times d_2$. By analogy with the one-dimensional case, we put $\mathbf{D}\vec{\varkappa} = {\rm cov}(\vec{\varkappa}, \vec{\varkappa})$.

Recall that $s \ge \underline{s}$, $N_{lb} | N$, and ${L_{sb}} | h$. Let
\begin{gather}
Z_i= f^{(s)}(\varepsilon_{h(i-1)+1}, \ldots, \varepsilon_{h(i-1)+s})
= f^{({\underline{s}})}(\varepsilon_{h(i-1)+1}, \ldots, \varepsilon_{h(i-1)+{\underline{s}}}), \quad i \ge 1, \nonumber \\
\Phi = \mathbf{D}_{H_0}Z_1 + \sum_{2 \le i \le \frac{s-1}h+1}\Big( \mathrm{cov}_{H_0}(Z_1, Z_i)
+ \big(\mathrm{cov}_{H_0}(Z_1, Z_i)\big)^{\mathrm{T}} \Big). \label{defPhi} \end{gather}
The matrix $\Phi$ is well defined. Indeed, we assume that the statistics $T_{sum}, T_{lb}$, and $T_{sb}$ are well defined, so the numbers $\sigma_{sum}, \sigma_{lb}, E_{sb}(0), \ldots, E_{sb}(K_{sb})$ are finite and nonzero. Since $\sigma_{sum} \in (0,\infty)$, the first term on the right-hand side of the formula \eqref{dlyaSigmaSquared} cannot be infinite, i.e. $\mathbf{D}_{H_0} f_{sum}(\varepsilon_1, \ldots, \varepsilon_{{m_{sum}}}) < \infty$. Similarly, $\mathbf{D}_{H_0} f_{lb}(\varepsilon_1, \ldots, \varepsilon_{{m_{lb}}}) < \infty$. As a consequence, the matrix $\mathbf{D}_{H_0}Z_1$ is well defined, and therefore so is the matrix $\Phi$.

Note that the quantities $N,h,s$ and $Z_i$ are used in the proof of the main results as follows: $N$ is the number of disjoint blocks into which the tested sequence $(\varepsilon_1,\ldots,\varepsilon_n)$ is splitted, $s$ is the length of (generally speaking, intersecting) chains of consecutive symbols into which each of the given $N$ blocks is splitted, with only each $h$-th of such $s$-chains being taken into account; the quantities $Z_i$ contain information about what “happens” in the $i$-th of such $s$-chains, i.e. in $(\varepsilon_{h(i-1)+1}, \ldots, \varepsilon_{h(i-1)+s})$.

For $1 \le k \le N$, $0 \le j \le K_{sb}+2$ put
\begin{gather}
\mathcal{X}_k^{(N,s,h)}(j) = \frac{1}{\sqrt{n}} \sum_{\frac{L(k-1)}h\le i \le \frac{Lk-s}h} \Big( f_j^{(s)}(\varepsilon_{hi+1}, \ldots, \varepsilon_{hi+s}) - \mathbf{E}_{H_0}f_j^{(s)}(\varepsilon_{1}, \ldots, \varepsilon_{s}) \Big), \nonumber \\
\mathcal{X}_k^{(N,s,h)} = \Big( \mathcal{X}_k^{(N,s,h)}(0), \mathcal{X}_k^{(N,s,h)}(1), \ldots, \mathcal{X}_k^{(N,s,h)}(K_{sb}+2) \Big). \nonumber
\end{gather}
As a corollary,
\begin{equation}
\mathcal{X}_k^{(N,s,h)} = \frac{1}{\sqrt{n}} \sum_{\frac{L(k-1)}h\le i \le \frac{Lk-s}h} ( Z_{i+1} - \mathbf{E}_{H_0}Z_1). \nonumber
\end{equation}
Put
\begin{gather}
\mathcal{Y}_k^{(N,s,h)}(j) = \sum_{i = 1}^{\frac{N}{N_{lb}}} \mathcal{X}_{\frac{N}{N_{lb}}(k-1)+i}^{(N,s,h)}(j), \quad 1 \le k \le N_{lb}, \ 0 \le j \le K_{sb}+2,\nonumber \\
\mathcal{Y}_k^{(N,s,h)} = \sum_{i = 1}^{\frac{N}{N_{lb}}} \mathcal{X}_{\frac{N}{N_{lb}}(k-1)+i}^{(N,s,h)}, \quad 1 \le k \le N_{lb}.\label{defYk}
\end{gather}

As a simple alternative hypothesis $H_1$ hereafter we will consider the following series scheme: the $n$-th series contains $n$ random variables $\varepsilon_1, \ldots, \varepsilon_n$, taking values in the set $\mathfrak{X}$ and having an unambiguously defined joint distribution.

\begin{remark}\label{zamechKlonH0}
The simple hypothesis $H_0$ can be treated as a special case of $H_1$. Indeed, the special case of hypothesis $H_1$ is such a scheme of series $H_{0,n}$ that the $n$-th series contains $n$ independent random variables $\varepsilon_1, \ldots, \varepsilon_n$ with distribution $\mathbf{P}_0$. Both $H_0$ and $H_{0,n}$ will be denoted by $H_0$ for simplicity.
\end{remark}

Next, for $1 \le k \le N$ let a vector
\begin{gather*}
\eta_k = \big(\eta_k(0), \eta_k(1), \ldots, \eta_k(K_{sb}+2) \big)
\end{gather*}
be a $(K_{sb}+3)$-dimensional random vector. We need the following condition on the alternative $H_1$: for fixed $N, s, h$, the following convergence holds:
\begin{gather}\Big(\mathcal{X}_1^{(N,s,h)}(0), \ldots, \mathcal{X}_1^{(N,s,h)}(K_{sb}+2), \ldots, \mathcal{X}_N^{(N,s,h)}(0), \ldots, \mathcal{X}_N^{(N,s,h)}(K_{sb}+2)\Big) \xrightarrow{d} \nonumber \\
\xrightarrow{d} \Big(\eta_1(0), \ldots, \eta_1(K_{sb}+2), \ldots, \eta_N(0), \ldots, \eta_N(K_{sb}+2)\Big), \quad n \to \infty. \nonumber
\end{gather}
Understanding the convergence of random vectors whose components are themselves random vectors in the same way as in \cite[Ch. 1, \S 4]{Billingsley}, the last relation can be written in a simplified form:
\begin{gather}\Big(\mathcal{X}_1^{(N,s,h)}, \ldots, \mathcal{X}_N^{(N,s,h)}\Big) \xrightarrow{d} (\eta_1, \ldots, \eta_N), \quad n \to \infty. \label{XshoditsyaKeta}
\end{gather}

\begin{example} \label{zamechPriH0Horosho1}
If $H_0$ is true, then the vectors $\mathcal{X}_1^{(N,s,h)}, \ldots, \mathcal{X}_N^{(N,s,h)}$ are independent (as functions of independent random variables $\varepsilon_i$), and for $k \in \{1,\ldots, N\}$ from the central limit theorem applied to $\big[\frac{s-1}h\big]$-dependent random vectors $Z_{i+1}$ \Big($ \frac{L(k-1)}h\le i \le \frac{Lk-s}h$\Big) it follows that the limit distribution of $\mathcal{X}_k^{(N,s,h)}$ is $\mathcal{N}(0, \frac{1}{Nh}\Phi)$, where the matrix $\Phi$ is defined in \eqref{defPhi}. Thus, if $H_0$ is true, then the relation \eqref{XshoditsyaKeta} holds for independent random vectors $\eta_1, \ldots, \eta_N$ such that $\eta_k \sim \mathcal{N}(0, \frac{1}{Nh}\Phi)$ (see also Example \ref{zamechPriH0Horosho3} below).
\end{example}

Next, put
\begin{gather}
\zeta_k(j) = \sum_{i = 1}^{\frac{N}{N_{lb}}} \eta_{\frac{N}{N_{lb}}(k-1)+i}(j), \quad 1 \le k \le N_{lb}, 0 \le j \le K_{sb}+2, \nonumber\\
\zeta_k = \sum_{i = 1}^{\frac{N}{N_{lb}}} \eta_{\frac{N}{N_{lb}}(k-1)+i}, \quad 1 \le k \le N_{lb}.\label{defZetaK}
\end{gather}

\begin{theorem} \label{th1}
Let $s \ge {\underline{s}}$, $N_{lb} | N$ and ${L_{sb}} | h$. Assume the numbers $N, s, h, m_{sum}, m_{lb}, N_{lb}, L_{sb}, K_{sb}$, the sets $\alpha_{sb}(0), \ldots, \alpha_{sb}$, the function $f_{sb}$ and the bounded functions $f_{sum}$ and $f_{lb}$ are fixed. Then
\begin{gather}
T_{sum} = \sum_{k=1}^{N_{lb}} \mathcal{Y}_k^{(N,s,h)}({K_{sb}}+2) +o_P(1), \quad n \to \infty. \label{flaTsumAsRazl}
\end{gather}
If, moreover, the relation \eqref{XshoditsyaKeta} is fulfilled, then
\begin{gather}
T_{lb} = N_{lb} \cdot \sum_{k=1}^{N_{lb}}  \Big(\mathcal{Y}_k^{(N,s,h)}({K_{sb}}+1) \Big)^2    +o_P(1), \label{flaTlbAsRazl} \\
T_{sb} = L_{sb} \cdot \sum_{j=0}^{{K_{sb}}} \bigg( \sum_{k=1}^{N_{lb}}  \mathcal{Y}_k^{(N,s,h)}(j) \bigg)^2  +o_P(1), \label{flaTsbAsRazl} \\
\Big(\mathcal{Y}_1^{(N,s,h)}, \ldots, \mathcal{Y}_{N_{lb}}^{(N,s,h)}\Big) \xrightarrow{d} (\zeta_1, \ldots, \zeta_{N_{lb}}) \label{limY}
\end{gather}
at $n \to \infty$, and, as a corollary,
\begin{gather}
\Big( T_{sum} ,  T_{lb},   T_{sb}  \Big) \xrightarrow{d} \nonumber \\
 \xrightarrow{d}  \Bigg( \sum_{k=1}^{N_{lb}} \zeta_{k}({K_{sb}}+2), N_{lb} \cdot \sum_{k=1}^{N_{lb}} \zeta_k^2({K_{sb}}+1),  L_{sb} \cdot \sum_{j=0}^{{K_{sb}}} \Bigg(\sum_{k=1}^{N_{lb}} \zeta_k(j) \Bigg)^2 \Bigg)\label{th1LimDistr}
\end{gather}
at $n \to \infty$.
\end{theorem}

\begin{corollary} \label{corolOsovmRaspr3stat}
Suppose that $s \ge {\underline{s}}$, $N_{lb} | N$, ${L_{sb}} | h$ and the hypothesis $H_0$ holds. Let random vectors $\eta_1, \ldots, \eta_N$ be independent and $\eta_k \sim \mathcal{N}\Big(0, \frac{1}{Nh}\Phi\Big)$, $1 \le k \le N$. Define vectors $\zeta_1, \ldots, \zeta_{N_{lb}}$ using the formula \eqref{defZetaK}. Assume that the numbers $N, s, h, m_{sum}, m_{lb}, N_{lb}, L_{sb}, K_{sb}$, the functions $f_{sb}, f_{sum}, f_{lb}$ and the sets $\alpha_{sb}(ts0), \ \alpha_{sb}(K_{sb})$ are fixed. Then the vectors $\zeta_1, \ldots, \zeta_{N_{lb}}$ are independent, $\zeta_k \sim \mathcal{N}\Big(0, \frac{1}{hN_{lb}}\Phi\Big)$ for $1 \le k \le N_{lb}$ and relations \eqref{flaTsumAsRazl}--\eqref{th1LimDistr} hold.
\end{corollary}

If the conditions of Corollary \ref{corolOsovmRaspr3stat} are satisfied, then by \eqref{defZetaK} and \eqref{th1LimDistr} the limit joint distribution of the statistics $T_{sum}, T_{lb}, T_{sb}$ coincides with the distribution of a known function (a vector of quadratic and linear forms) of independent random vectors $\eta_k, 1 \le k \le N,$ with known distributions: $\eta_k \sim \mathcal{N}\Big(0, \frac{1}{Nh}\Phi\Big), 1 \le k \le N$.

Next, we consider the question of the limit joint distribution of several statistics of different types. Consider $Q$ triples of statistics $\Big(T_{sum}^{[q]}, T_{lb}^{[q]}, T_{sb}^{[q]}\Big), 1 \le q \le Q$, where $T_{sum}^{[q]}, T_{lb}^{[q]}$, and $T_{sb}^{[q]}$ are the summing, long-block, and short-block statistics, respectively. We will use the following notation.
All parameters related to $T_{sum}^{[q]}$ will be marked with the subscript ''sum'' and the superscript $[q]$:
$$T_{sum}^{[q]} = \frac{ \sum_{i=1}^{n-m_{sum}^{[q]}+1} \Big(f_{sum}^{[q]}\big(\varepsilon_i, \ldots, \varepsilon_{i+m_{sum}^{[q]}-1}\big) - E_{sum}^{[q]}\Big)}{\sigma_{sum}^{[q]} \sqrt{n - m_{sum}^{[q]}+1}}.$$
Similarly, $N_{lb}^{[q]}$ is the number of long blocks used to construct the statistics $T_{lb}^{[q]}$, $L_{sb}^{[q]}$ --- the length of short blocks used in constructing the statistic $T_{sb}^{[q]}$, etc. The triple of statistics $\Big(T_{sum}^{[q]}, T_{lb}^{[q]}, T_{sb}^{[q]}\Big)$ is constructed from the set of functions $(f_{sum}^{[q]} f_{lb}^{[q]} f_{sb}^{[q]})$. Quantities whose values depend on the parameters of the statistics $T_{sum}^{[1]}, T_{lb}^{[1]}, T_{sb}^{[1]}, \ldots, T_{sum}^{[Q]}, T_{lb}^{[Q]}, T_{sb}^{[Q]}$ will be marked with an asterisk.

Put
$$ \underline{s}^* = h + \max_{1 \le q \le Q} \max\Big(m_{sum}^{[q]}, m_{lb}^{[q]}\Big) -1.$$
Assume that $N_{lb}^{[q]} | N$ and $L_{sb}^{[q]} | h$ for all $1 \le q \le Q$.

Let $$K^* = 3Q+\sum_{q=1}^Q K_{sb}^{[q]}.$$
We define the function $f^{*(\underline{s}^*)}: \mathfrak{X}^{\underline{s}^*} \to \mathbb{R}^{K^*}$. To do this, we define the quantities $f_j^{*(\underline{s}^*)}, 0 \le j \le K^*-1$, which are the coordinates of the vector $f^{*(\underline{s}^*)}$. Put
\begin{gather}
f_j^{*({\underline{s}^*})}\big(\theta_{1}, \ldots, \theta_{{\underline{s}^*}}\big)
=\begin{cases}
\frac{1}{\sqrt{E_{sb}^{[1]}(j)}}\sum_{u=1}^{h / L_{sb}^{[1]}}I_{f_{sb}^{[1]}\big(\theta_{L_{sb}^{[1]}(u-1)+1}, \ldots, \theta_{{L_{sb}^{[1]}}u}\big)\in \alpha_{sb}^{[1]}(j) } ,&\text{if $0 \le j \le {K_{sb}^{[1]}}$,}\\
\frac{1}{\sigma_{lb}^{[1]}} \cdot \sum_{u=1}^{h} f_{lb}^{[1]}\big(\theta_{u}, \ldots, \theta_{u+{m_{lb}^{[1]}}-1}\big),&\text{if $j = K_{sb}^{[1]}+1$,}\\
\frac{1}{\sigma_{sum}^{[1]}} \cdot \sum_{u=1}^{h} f_{sum}^{[1]}\big(\theta_u, \ldots, \theta_{u+m_{sum}^{[1]}-1}\big),&\text{if $j = K_{sb}^{[1]}+2$.}
\end{cases} \nonumber
\end{gather}
Let
\begin{gather}
f_{K_{sb}^{[1]}+3+j}^{*({\underline{s}^*})}\big(\theta_{1}, \ldots, \theta_{{\underline{s}^*}}\big)
=\begin{cases}
\frac{1}{\sqrt{E_{sb}^{[2]}(j)}}\sum_{u=1}^{h / L_{sb}^{[2]}}I_{f_{sb}^{[2]}\big(\theta_{L_{sb}^{[2]}(u-1)+1}, \ldots, \theta_{{L_{sb}^{[2]}}u}\big)\in \alpha_{sb}^{[2]}(j) } ,&\text{if $0 \le j \le {K_{sb}^{[2]}}$,}\\
\frac{1}{\sigma_{lb}^{[2]}} \cdot \sum_{u=1}^{h} f_{lb}^{[2]}\big(\theta_{u}, \ldots, \theta_{u+{m_{lb}^{[2]}}-1}\big),&\text{if $j = K_{sb}^{[2]}+1$,}\\
\frac{1}{\sigma_{sum}^{[2]}} \cdot \sum_{u=1}^{h} f_{sum}^{[2]}\big(\theta_u, \ldots, \theta_{u+m_{sum}^{[2]}-1}\big),&\text{if $j = K_{sb}^{[2]}+2$.}
\end{cases} \nonumber
\end{gather}
Similarly, we define $f_{K_{sb}^{[1]}+K_{sb}^{[2]}+6+j}^{*({\underline{s}^*})}\big(\theta_{1}, \ldots, \theta_{{\underline{s}^*}}\big)$ for $0 \le j \le K_{sb}^{[3]}+2$, then we define $f_{K_{sb}^{[1]}+K_{sb}^{[2]}+K_{sb}^{[3]}+9+j}^{*({\underline{s}^*})}\big(\theta_{1}, \ldots, \theta_{{\underline{s}^*}}\big)$ for $0 \le j \le K_{sb}^{[4]}+2$, etc. Note that the first $K_{sb}^{[1]}+3$ coordinates of the vector $f^{*(\underline{s}^*) }$ correspond to statistics
$T_{sum}^{[1]}, T_{lb}^{[1]}, T_{sb}^{[1]}$, the next $K_{sb}^{[2]}+3$ coordinates correspond to statistics $T_{sum}^{[2]}, T_{lb}^{[2]}, T_{sb}^{[2]}$, etc., the last $K_{sb}^{[Q]}+3$ coordinates correspond to statistics $T_{sum}^{[Q]}, T_{lb}^{[Q]}, T_{sb}^{[Q]}$.

Next, let $s \ge {\underline{s}^*}$. For each $(\theta_{1}, \ldots, \theta_{s}) \in \mathfrak{X}^s$ put
\begin{gather}
f_j^{*(s)}(\theta_{1}, \ldots, \theta_{s}) = f_j^{*({\underline{s}^*})}(\theta_{1}, \ldots, \theta_{{\underline{s}^*}}), \quad 0 \le j \le K^*-1, \nonumber \\
f^{*(s)}(\theta_{1}, \ldots, \theta_{s}) = f^{*({\underline{s}^*})}\big(\theta_{1}, \ldots, \theta_{{\underline{s}^*}}\big).\nonumber
\end{gather}
Recall that $N_{lb}^{[q]} | N$ and $L_{sb}^{[q]} | h$ for all $1 \le q \le Q$ and the inequality $s \ge \underline{s}^{*}$ holds. Let
\begin{gather}
Z_i^*= f^{*(s)}(\varepsilon_{h(i-1)+1}, \ldots, \varepsilon_{h(i-1)+s})
= f^{*({\underline{s}^*})}(\varepsilon_{h(i-1)+1}, \ldots, \varepsilon_{h(i-1)+{\underline{s}^*}}), \quad i \ge 1, \nonumber \\
\Phi^* = \mathbf{D}_{H_0}Z_1^* + \sum_{2 \le i \le \frac{s-1}h + 1}\Big( \mathrm{cov}_{H_0}(Z_1^*, Z_i^*)
+ \big(\mathrm{cov}_{H_0}(Z_1^*, Z_i^*)\big)^{\top} \Big). \nonumber \end{gather}
The matrix $\Phi^*$ is well defined by reasons similar to those used to prove that the matrix $\Phi$ is well defined.

For $1 \le k \le N$, $0 \le j \le K^*-1$ put
\begin{gather}
\mathcal{X}_k^{*(N,s,h)}(j) = n^{-\frac12}\!\!\!\!\!\!\!\sum_{\frac{L(k-1)}h\le i \le \frac{Lk-s}h} \!\!\!\!\Big( f_j^{*(s)}(\varepsilon_{hi+1}, \ldots, \varepsilon_{hi+s}) - \mathbf{E}_{H_0}f_j^{*(s)}(\varepsilon_1,\ldots, \varepsilon_s)\Big), \nonumber \\
\mathcal{X}_k^{*(N,s,h)} = \bigg( \mathcal{X}_k^{*(N,s,h)}(0), \mathcal{X}_k^{*(N,s,h)}(1), \ldots, \mathcal{X}_k^{*(N,s,h)}(K^*-1)\bigg). \nonumber
\end{gather}
As a corollary,
\begin{equation}
\mathcal{X}_k^{*(N,s,h)} = \frac{1}{\sqrt{n}} \sum_{\frac{L(k-1)}h\le i \le \frac{Lk-s}h} ( Z_{i+1}^* - \mathbf{E}_{H_0}Z_1^* ). \nonumber
\end{equation}
Let
\begin{gather}
\mathcal{Y}_k^{[1],(N,s,h)}(j) = \sum_{i = 1}^{\frac{N}{N_{lb}^{[1]}}} \mathcal{X}_{\frac{N}{N_{lb}^{[1]}}(k-1)+i}^{*(N,s,h)}(j),\nonumber \\
\mathcal{Y}_k^{[1],(N,s,h)} = \Big( \mathcal{Y}_k^{[1],(N,s,h)}(0), \mathcal{Y}_k^{[1],(N,s,h)}(1), \ldots, \mathcal{Y}_k^{[1],(N,s,h)}\Big(K_{sb}^{[1]}+2\Big)\Big) \nonumber
\end{gather}
for $1 \le k \le N_{lb}^{[1]}, \ 0 \le j \le K_{sb}^{[1]}+2$.

Put
\begin{gather}
\mathcal{Y}_k^{[2],(N,s,h)}(j) = \sum_{i = 1}^{\frac{N}{N_{lb}^{[2]}}} \mathcal{X}_{\frac{N}{N_{lb}^{[2]}}(k-1)+i}^{*(N,s,h)}\Big(K_{sb}^{[1]}+3+j\Big), \nonumber \\
\mathcal{Y}_k^{[2],(N,s,h)} = \Big( \mathcal{Y}_k^{[2],(N,s,h)}(0), \mathcal{Y}_k^{[2],(N,s,h)}(1), \ldots, \mathcal{Y}_k^{[2],(N,s,h)}\Big(K_{sb}^{[2]}+2\Big)\Big) \nonumber
\end{gather}
for $1 \le k \le N_{lb}^{[2]}, \ 0 \le j \le K_{sb}^{[2]}+2$. Similarly, for $q \in \{3, \ldots, Q\}$ we define the quantities $\mathcal{Y}_k^{[q],(N,s,h)}(j)$ and $\mathcal{Y}_k^{[q],(N,s,h)}$ for $1 \le k \le N_{lb}^{[q]}, \ 0 \le j \le K_{sb}^{[q]}+2$.

Note that in \cite{Savelov11} the quantities $X_k^{*(N,s,h)}$ and $Y_k^{[q],(N,s,h)}$ were defined, similar to the quantities $ \mathcal{X}_k^{*(N,s,h)}$ and $\mathcal{Y}_k^{[q],(N,s,h)}$ from the present paper.

Next, for $1 \le k \le N$, we consider $K^*$-dimensional random vectors
\begin{equation}
\eta_k^* = \bigg(\eta_k^*(0), \eta_k^*(1), \ldots, \eta_k^*(K^*-1) \bigg).\nonumber
\end{equation}
We need the following condition on the alternative $H_1$: for fixed $N, s, h$, the convergence holds
\begin{equation}\Big(\mathcal{X}_1^{*(N,s,h)}, \ldots, \mathcal{X}_N^{*(N,s,h)}\Big) \xrightarrow{d} (\eta_1^*, \ldots, \eta_N^*), \quad n \to \infty. \label{XStarshoditsyaKetaStar}
\end{equation}

\begin{example} \label{zamechPriH0Horosho2}
Arguing in the same way as in Example \ref{zamechPriH0Horosho1}, we find that if $H_0$ is true, then the relation \eqref{XStarshoditsyaKetaStar} is satisfied for independent random vectors $\eta_1^*, \ldots, \eta_N^*$ such that $\eta_k^* \sim \mathcal{N}(0, \frac{1}{Nh}\Phi^*)$. In other words, $H_0$ is a special case of the alternative $H_1$ satisfying the condition \eqref{XStarshoditsyaKetaStar} (see Remark \ref{zamechKlonH0}).
\end{example}

\begin{example} \label{zamechPriH0Horosho3}
If the elements of the tested sequence $\varepsilon_1, \ldots, \varepsilon_n$ take values in a fixed finite set $\{0,1,\ldots,R-1\}$ and all the conditions of Theorem 2 \cite{Savelov11} are satisfied (they guarantee that $H_1$ is in some sense close to $H_0$), then the relation \eqref{XStarshoditsyaKetaStar} is satisfied for the corresponding vectors $\eta_1^*, \ldots, \eta_N^*$ from Theorem 2 \cite{Savelov11}, since in this case the relation \cite[(28)]{Savelov11} holds and for any $k \in \{1,\ldots,N\}$ the vector $\mathcal{X}_k^{*(N,s,h)}$ coincides with the vector $X_k^{*(N,s,h)}$ from \cite{Savelov11} up to $o_P(1)$. We also note that examples of Markov alternatives $H_1$ for which all the conditions of Theorem 2 \cite{Savelov11} are satisfied are considered in \cite{Savelov11} (see hypotheses $\tilde{H}_1$ and $\hat{H}_1$).
\end{example}

\begin{remark} \label{zamechProH1BlizkH0}
Examples \ref{zamechPriH0Horosho2} and \ref{zamechPriH0Horosho3} show that the relation \eqref{XStarshoditsyaKetaStar} can be interpreted as a condition meaning that $H_1$ and $H_0$ are close in some sense.
\end{remark}

Recall that we assume that $N_{lb}^{[q]} | N$ and $L_{sb}^{[q]} | h$ for all $1 \le q \le Q$ and the inequality $s \ge \underline{s}^{*}$ holds. Put
\begin{gather}
\zeta_k^{[1]}(j) = \sum_{i = 1}^{\frac{N}{N_{lb}^{[1]}}} \eta_{\frac{N}{N_{lb}^{[1]}}(k-1)+i}^*(j), \nonumber \\
\zeta_k^{[1]} = \bigg( \zeta_k^{[1]}(0), \zeta_k^{[1]}(1), \ldots, \zeta_k^{[1]}\Big(K_{sb}^{[1]}+2\Big)\bigg) \nonumber
\end{gather}
for $1 \le k \le N_{lb}^{[1]}, \ 0 \le j \le K_{sb}^{[1]}+2$.

Let
\begin{gather}
\zeta_k^{[2]}(j) = \sum_{i = 1}^{\frac{N}{N_{lb}^{[2]}}} \eta_{\frac{N}{N_{lb}^{[2]}}(k-1)+i}^*\Big(K_{sb}^{[1]}+3+j\Big),\nonumber \\
\zeta_k^{[2]} = \Big( \zeta_k^{[2]}(0), \zeta_k^{[2]}(1), \ldots, \zeta_k^{[2]}\Big(K_{sb}^{[2]}+2\Big)\Big)\nonumber
\end{gather}
for $1 \le k \le N_{lb}^{[2]}, \ 0 \le j \le K_{sb}^{[2]}+2$. Similarly, for $q \in \{3, \ldots, Q\}$, we define the quantities $\zeta_k^{[q]}(j)$ and $\zeta_k^{[q]}$ for $1 \le k \le N_{lb}^{[q]}, \ 0 \le j \le K_{sb}^{[q]}+2$.

Consider $J$ quadratic statistics $T_{quad}^{[1]}, \ldots, T_{quad}^{[J]}$, constructed from $T_{sum}^{[1]}, \ldots, T_{sum}^{[Q]}$:
\begin{equation}
T_{quad}^{[j]} = \sum_{i=1}^{\tau_{quad}^{[j]}} \Bigg(\sum_{q=1}^{Q} d_{quad}^{[j]}(i,q) T_{sum}^{[q]}\Bigg)^2. \label{deTquadJ}
\end{equation}

Recall that by the parameters of the statistics $T_{sum}^{[q]}, T_{lb}^{[q]}, T_{sb}^{[q]}$ we mean the numbers $m_{sum}^{[q]}, m_{lb}^{[q]}, N_{lb}^{[q]}, L_{sb}^{[q]}, K_{sb}^{[q]}$, the functions $f_{sum}^{[q]}, f_{lb}^{[q]}, f_{sb}^{[q]}$, and the sets $\alpha_{sb}^{[q]}(0), \ldots, \alpha_{sb}^{[q]}\Big(K_{sb}^{[q]}\Big)$. In this case, the quantities $L_{lb}^{[q]}$ and $N_{sb}^{[q]}$ are expressed in terms of the parameters of the statistics and the value $n$, but are not called parameters themselves. By the parameters of the statistic $T_{quad}^{[j]}$ we mean the numbers $Q, \tau_{quad}^{[j]}$ and $d_{quad}^{[j]}(i,q)$ ($1 \le i \le \tau_{quad}^{[j]}, 1 \le q \le Q$), as well as the parameters of the statistics $T_{sum}^{[1]}, \ldots, T_{sum}^{[Q]}$.

\begin{theorem} \label{th2}
Let $s \ge \underline{s}^*$ and let the relations $N_{lb}^{[q]} | N$ and $L_{sb}^{[q]} | h$ hold for all $q \in \{1, \ldots, Q\}$. Suppose that the numbers $N, s, h$ and all parameters of the statistics $T_{sum}^{[1]}, T_{lb}^{[1]}, T_{sb}^{[1]}, \ldots, T_{sum}^{[Q]}, T_{lb}^{[Q]}, T_{sb}^{[Q]}, T_{quad}^{[1]},\ldots, T_{quad}^{[J]}$ are fixed. Let the functions $f_{sum}^{[q]}$ and $f_{lb}^{[q]}$ be bounded for $1 \le q \le Q$. Then
\begin{gather}
T_{sum}^{[q]} = \sum_{k=1}^{N_{lb}^{[q]}} \mathcal{Y}_k^{[q],(N,s,h)}\Big({K_{sb}^{[q]}}+2\Big) +o_P(1), \quad n \to \infty. \label{flaTsumAsRazl2}
\end{gather}
If, in addition, the relation \eqref{XStarshoditsyaKetaStar} holds, then
\begin{gather}
T_{lb}^{[q]} = N_{lb}^{[q]} \cdot \sum_{k=1}^{N_{lb}^{[q]}} \Big(\mathcal{Y}_k^{[q],(N,s,h)}\Big(K_{sb}^{[q]}+1\Big) \Big)^2 +o_P(1), \label{flaTlbAsRazl2} \\
T_{sb}^{[q]} = L_{sb}^{[q]} \cdot \sum_{j=0}^{K_{sb}^{[q]}} \bigg( \sum_{k=1}^{N_{lb}^{[q]}} \mathcal{Y}_k^{[q],(N,s,h)}(j) \bigg)^2 +o_P(1), \label{flaTsbAsRazl2} \\
T_{quad}^{[j]} = \sum_{i=1}^{\tau_{quad}^{[j]}} \Bigg(\sum_{q=1}^{Q} d_{quad}^{[j]}(i,q) \sum_{k=1}^{N_{lb}^{[q]}} \mathcal{Y}_k^{[q],(N,s,h)}\Big({K_{sb}^{[q]}}+2\Big) \Bigg)^2 +o_P(1),
\label{flaTquadAsRazl}\\
\Big(\mathcal{Y}_1^{[1],(N,s,h)}, \ldots, \mathcal{Y}_{N_{lb}^{[1]}}^{[1],(N,s,h)}, \ldots, \mathcal{Y}_1^{[Q],(N,s,h)}, \ldots, \mathcal{Y}_{N_{lb}^{[Q]}}^{[Q],(N,s,h)}\Big) \xrightarrow{d}
\nonumber\\
\xrightarrow{d}\Big(\zeta_1^{[1]}, \ldots, \zeta_{N_{lb}^{[1]}}^{[1]}, \ldots, \zeta_1^{[Q]}, \ldots, \zeta_{N_{lb}^{[Q]}}^{[Q]}\Big) \label{limYbig}
\end{gather}
as $n \to \infty$, which implies that
\begin{gather}
\Big( T_{sum}^{[1]} , T_{lb}^{[1]}, T_{sb}^{[1]}, \ldots, T_{sum}^{[Q]} , T_{lb}^{[Q]}, T_{sb}^{[Q]}, T_{quad}^{[1]},\ldots, T_{quad}^{[J]} \Big) \xrightarrow{d} \nonumber \\
\xrightarrow{d} \Bigg( \sum_{k=1}^{N_{lb}^{[1]}} \zeta_{k}^{[1]}\Big({K_{sb}^{[1]}}+2\Big), N_{lb}^{[1]} \sum_{k=1}^{N_{lb}^{[1]}} \Big(\zeta_k^{[1]}\Big(K_{sb}^{[1]}+1\Big)\Big)^2, L_{sb}^{[1]} \sum_{j=0}^{K_{sb}^{[1]}} \Bigg(\sum_{k=1}^{N_{lb}^{[1]}} \zeta_k^{[1]}(j) \Bigg)^2, \ldots \nonumber \\
\ldots, \sum_{k=1}^{N_{lb}^{[Q]}} \zeta_{k}^{[Q]}\Big(K_{sb}^{[Q]}+2\Big), N_{lb}^{[Q]} \sum_{k=1}^{N_{lb}^{[Q]}} \Big(\zeta_k^{[Q]}\Big(K_{sb}^{[Q]}+1\Big)\Big)^2, L_{sb}^{[Q]} \sum_{j=0}^{K_{sb}^{[Q]}} \Bigg(\sum_{k=1}^{N_{lb}^{[Q]}} \zeta_k^{[Q]}(j) \Bigg)^2, \nonumber \\
\theta^{[1]}, \ldots, \theta^{[J]}\Bigg)\label{th2LimDistrGeneral}
\end{gather}
at $n \to \infty$, where
\begin{equation*}
\theta^{[j]} = \sum_{i=1}^{\tau_{quad}^{[j]}} \Bigg(\sum_{q=1}^{Q} d_{quad}^{[j]}(i,q) \sum_{k=1}^{N_{lb}^{[q]}} \zeta_{k}^{[q]}\Big({K_{sb}^{[q]}}+2\Big) \Bigg)^2, \quad 1 \le j \le J.
\end{equation*}
\end{theorem}
According to Theorem \ref{th2} random vectors $\mathcal{X}_1^{*(N,s,h)}, \ldots, \mathcal{X}_N^{*(N,s,h)}$ are a kind of "basis" in the sense that each of the statistics $T_{sum}^{[q]}, T_{lb}^{[q]}, T_{sb}^{[q]}, T_{quad}^{[j]}$ ($1 \le q \le Q, 1 \le j \le J$) can be represented as a function (namely, as a linear or quadratic form) of these "basis" vectors with error $o_P(1)$, $n \to \infty$. This representation is used in the proof of Theorem \ref{th2} to prove the formula \eqref{th2LimDistrGeneral}, i.e. to find the limit joint distribution of summing, long-block, short-block, and quadratic statistics. The statement of Theorem \ref{th2} overlaps with the results of \cite{Savelov11, Savelov12}.

Note that Theorem \ref{th2} is in some sense more general than Theorem 2 \cite{Savelov11}, since the latter is applicable only to discrete test sequences. Moreover, Theorem 2 \cite{Savelov11} proves relation \cite[(28)]{Savelov11}, which is analogous (see Example \ref{zamechPriH0Horosho3}) to the relation \eqref{XStarshoditsyaKetaStar} from this paper, which appears in the conditions, and not in the conclusion of Theorem \ref{th2}. Roughly speaking, Theorem \ref{th2} is more general due to the fact that the relation \eqref{XStarshoditsyaKetaStar} is not proved, but acts as one of its conditions.

The statement of Theorem \ref{th2} involves $Q$ triples of statistics $\Big(T_{sum}^{[q]}, T_{lb}^{[q]}, T_{sb}^{[q]}\Big)$, i.e., the numbers of summing, long-block, and short-block statistics are equal. The more general case, when the numbers of statistics of these three types are not necessarily equal, can be easily reduced to the previous one. To do this, it is sufficient to supplement the set of statistics under consideration, apply Theorem \ref{th2}, and then "discard" the redundant statistics.

\begin{remark} \label{normalnoEsliJravnoNul}
As it was shown in the proof of Theorem \ref{th2}, Theorem \ref{th2} is also applicable when quadratic statistics are absent (i.e., $J=0$). Moreover, in the formulation of Theorem \ref{th2} for the case $J=0$, the vector \eqref{th2LimDistrGeneral} lacks components of the form $\theta^{[j]}$, and the equality \eqref{flaTquadAsRazl} is also absent. The statement of Theorem \ref{th1} is obtained (up to a change of notation) from Theorem \ref{th2} by substituting $Q=1$ and $J=0$.
\end{remark}

\begin{corollary} \label{thOsovmRaspr3stat}
Suppose that the hypothesis $H_0$ is true, $s \ge \underline{s}^*$, and for all $q \in \{1, \ldots, Q\}$ the relations $N_{lb}^{[q]} | N$ and $L_{sb}^{[q]} | h$ hold. Let random vectors $\eta_1^*, \ldots, \eta_N^*$ be independent and $\eta_k^* \sim \mathcal{N}(0, \frac{1}{Nh}\Phi^*)$, $1 \le k \le N$. Using $\eta_1^*, \ldots, \eta_N^*$, construct vectors $\zeta_1^{[q]}, \ldots, \zeta_{N_{lb}^{[q]}}^{[q]}$ for $1 \le q \le Q$ in the same way as above. Suppose that the numbers $N, s, h$ and all parameters of the statistics $T_{sum}^{[1]}, T_{lb}^{[1]}, T_{sb}^{[1]}, \ldots, T_{sum}^{[Q]}, T_{lb}^{[Q]}, T_{sb}^{[Q]}, T_{quad}^{[1]},\ldots, T_{quad}^{[J]}$ are fixed. Then the relations \eqref{flaTsumAsRazl2}--\eqref{th2LimDistrGeneral} hold.
\end{corollary}

Note that Corollary \ref{thOsovmRaspr3stat}, unlike Theorem \ref{th2}, does not assume that the functions $f_{sum}^{[q]}$ and $f_{lb}^{[q]}$ are bounded for $1 \le q \le Q$. A similar remark is also true for Corollary \ref{corolOsovmRaspr3stat} and Theorem \ref{th1}.

If the conditions of Corollary \ref{thOsovmRaspr3stat} are satisfied, then by virtue of \eqref{th2LimDistrGeneral} and the definition of the vectors $\zeta_1^{[q]}, \ldots, \zeta_{N_{lb}^{[q]}}^{[q]}$ ($1 \le q \le Q$), the limit joint distribution of the statistics $T_{sum}^{[1]} , T_{lb}^{[1]}, T_{sb}^{[1]}, \ldots, T_{sum}^{[Q]}, T_{lb}^{[Q]}, T_{sb}^{[Q]} $ coincides with the distribution of a known function (a vector of quadratic and linear forms) of independent random vectors $\eta_k^*, 1 \le k \le N,$ with known distributions: $\eta_k^* \sim \mathcal{N}\big(0, \frac{1}{Nh}\Phi^*\big), 1 \le k \le N$.

If the hypothesis $H_0^{U[0,1]}$ is true, Corollary \ref{thOsovmRaspr3stat} allows us to find the limit joint distribution of those statistics of the TestU01 package that are summing, long-block, short-block, or quadratic statistics (see Examples \ref{primerWeightDistrib}--\ref{primersmarsaSerialOver}). Similarly, if the hypothesis $H_0^{Bern(\frac12)}$ is true, Corollary \ref{thOsovmRaspr3stat} allows us to find the limit joint distribution of ten statistics of the following nine NIST STS tests: ''Monobit Test'', ''Frequency Test within a Block'', ''Runs Test'', ''Test for the Longest Run of Ones in a Block'', ''Binary Matrix Rank Test'', ''Non-overlapping Template Matching Test'', ''Linear Complexity Test'', ''Serial Test'' (two statistics correspond to this test), and ''Approximate Entropy Test'' (see a similar result in Remark 4 \cite{Savelov12}); each of these 10 statistics, as noted above, coincides up to $o_P(1)$ with a statistic of one of the four types.

Further, Theorem \ref{th1} allows us to obtain the limit distribution of the vector $\big(T_{sum} , T_{lb}, T_{sb} \big)$, and from the conditions of this theorem it follows that $\frac{1}{\sqrt{n}} \mathcal{X}_i^{(N,s,h)}(j) \overset{\mathbf{P}}{\to} 0$ for all $i \in \{1, \ldots ,N\}, \ j \in \{0, \ldots, K_{sb}+2\}$. In particular, the convergence of $\frac{1}{\sqrt{n}} \mathcal{X}_i^{(N,s,h)}(j) \overset{\mathbf{P}}{\to} 0$ holds if $H_0$ is true (see Example \ref{zamechPriH0Horosho1}). The case where this convergence does not hold is considered in the next lemma.

\begin{lemma} \label{LemmaProNeBlizkH1}
Let $s \ge {\underline{s}}$, $N_{lb} | N$, and ${L_{sb}} | h$. Assume that the numbers $N, s, h, m_{sum}, m_{lb}, N_{lb}, L_{sb}, K_{sb}$, the sets $\alpha_{sb}(0), \ldots, \alpha_{sb}(K_{sb})$, the function $f_{sb}$ and the bounded functions $f_{sum}$ and $f_{lb}$ are fixed. Let an alternative $H_1$ be true such that for all $1 \le i \le N$, $0 \le j \le K_{sb}+2$ the convergence $\frac{1}{\sqrt{n}} \mathcal{X}_i^{(N,s,h)}(j) \overset{\mathbf{P}}{\to} x_i(j)$, $n \to \infty$ holds, where $x_1(0), \ldots, x_N(K_{sb}+2)$ are some real numbers. For $1 \le k \le N_{lb}$, $0 \le j \le K_{sb}+2$ put
\begin{gather}
y_k(j) = \sum_{i = 1}^{\frac{N}{N_{lb}}}x_{\frac{N}{N_{lb}}(k-1)+i}(j), \nonumber \\
c_{sum} = \sum_{k=1}^{N_{lb}} y_k({K_{sb}}+2), \ \
c_{lb} = \sum_{k=1}^{N_{lb}} \big(y_k({K_{sb}}+1) \big)^2, \ \
c_{sb} = \sum_{j=0}^{K_{sb}} \bigg( \sum_{k=1}^{N_{lb}} y_k(j) \bigg)^2.\nonumber
\end{gather}
If $c_{sum} \ne 0$, then $T_{sum} \overset{\mathbf{P}}{\to} ({\rm{sgn}} (c_{sum}))\cdot(+\infty), \ n \to \infty$. If $c_{lb} \ne 0$, then $T_{lb }\overset{\mathbf{P}}{\to} +\infty,\ n \to \infty$. If $c_{sb} \ne 0$, then $T_{sb} \overset{\mathbf{P}}{\to} +\infty, \ n \to \infty$.
\end{lemma}

\begin{remark}
As Lemma \ref{LemmaProNeBlizkH1} shows (see also Remark \ref{zamechProH1BlizkH0}), if the alternative  $H_1$ is true, such that $H_1$ is ''not close'' to $H_0$, then in a wide class of cases the statistics $T_{sum}, T_{lb}, T_{sb}$ tend to infinity in absolute value. A similar result (albeit in a somewhat more cumbersome form) can easily be obtained for $T_{quad}$.
\end{remark}

If in the statement of Lemma \ref{LemmaProNeBlizkH1} we replace convergence in probability with almost sure convergence everywhere, the statement of the lemma remains valid (see Remark \ref{ProLemmyProNeblizkH1}). The statement of Lemma \ref{LemmaProNeBlizkH1} complements the statement of Lemma 5 \cite{Savelov11}.

\subsection{Analogue of the Berry-Esseen inequality}

Several results from the paper \cite{Savelov11}, which assumes that the elements of the tested sequence $\varepsilon_1, \ldots, \varepsilon_n$ take values in a fixed finite set $\{0,1,\ldots,R-1\}$ (i.e., the sample is discrete), were generalized above to the case where this condition is absent. For brevity, we will say that the results were generalized from the case $\mathfrak{X} = \{0,1,\ldots,R-1\}$ to the case of an arbitrary $\mathfrak{X}$. The results of the paper \cite{Savelov13}, in which an analogue of the Berry-Esseen inequality for summing, long-block, and short-block statistics was obtained (as well as the most important results of the paper \cite{Savelov12} --- see Subsection \ref{zamechOTomChtoVseOshchaetsya} below), are based on \cite{Savelov11} and are also transferred from the case of $\mathfrak{X} = \{0,1,\ldots,R-1\}$ to the case of an arbitrary $\mathfrak{X}$ with natural variations (we will take into account that in \cite{Savelov13} there is a typo: everywhere in the text, instead of the greatest common divisor, denoted by gcd, there should be the least common multiple, which we will denote by lcm). Namely, the statement of Theorem 2 \cite{Savelov13} (which assumes that $m_{sum}^{[q]} = m_{lb}^{[q]}=1$ for $1 \le q \le Q$) will be true for an arbitrary $\mathfrak{X}$ if the quantities
$M_{sum}^{[q]} = \max_{\theta \in \{0,1,\ldots,R-1\}}|f_{sum}^{[q]}(\theta)-E_{sum}^{[q]}|, M_{lb}^{[q]} = \max_{\theta \in \{0,1,\ldots,R-1\}}|f_{lb}^{[q]}(\theta)-E_{lb}^{[q]}|, N_{gcd}=gcd(N_{lb}^{[1]}, \ldots, N_{lb}^{[Q]}), h_{gcd}=gcd(L_{sb}^{[1]}, \ldots, L_{sb}^{[Q]})$
from \cite{Savelov13} are replaced by the quantities
$M_{sum}^{[q]} = \sup_{\theta \in \mathfrak{X}}|f_{sum}^{[q]}(\theta)-E_{sum}^{[q]}|, M_{lb}^{[q]} = \sup_{\theta \in \mathfrak{X}}|f_{lb}^{[q]}(\theta)-E_{lb}^{[q]}|, N_{lcm}=lcm(N_{lb}^{[1]}, \ldots, N_{lb}^{[Q]}), h_{lcm}=lcm(L_{sb}^{[1]}, \ldots, L_{sb}^{[Q]})$
(all other quantities: $M_{sb}^{[q]}(j), u_{sum}^{[q]}, \underline{n}, \delta_{BE}^*(N,h)$, etc. are defined in exactly the same way as in \cite{Savelov13}) and supplement the conditions of Theorem 2 \cite{Savelov13} with the condition of boundedness of the functions $f_{sum}^{[q]}$ and $f_{lb}^{[q]}$ for $1 \le q \le Q$ (which guarantees the finiteness of $M_{sum}^{[q]}$ and $M_{lb}^{[q]}$) and the condition $s=h$ (which is implicitly implied in \cite{Savelov13}, since in \cite{Savelov13} the function $f^{*(s)}$ and the quantities $\zeta_i^{[q]}$ are defined only for $s=h$). In other words, the following theorem holds.
\begin{theorem} \label{th2Berry}
Suppose that $s = h$ and let $N_{lb}^{[q]} | N, L_{sb}^{[q]} | h$, $m_{sum}^{[q]} = m_{lb}^{[q]}=1$ and $u_{sum}^{[q]}, u_{lb}^{[q]}, u_{sb}^{[q]} \in [0, +\infty]$ hold for all $q \in \{1, \ldots, Q\}$. Suppose that $f_{sum}^{[q]}$ and $f_{lb}^{[q]}$ are bounded for $1 \le q \le Q$, the inequality $n \ge N(2h-1)$ holds, and the hypothesis $H_0$ is true. Then the statistics $T_{sum}^{[q]}, T_{lb}^{[q]}, T_{sb}^{[q]}$ ($1 \le q \le Q$) are well defined and
\begin{gather}
\Bigg| \mathbf{P}_{H_0}\Big( \Big| T_{sum}^{[1]} \Big| \le u_{sum}^{[1]}, T_{lb}^{[1]} \le u_{lb}^{[1]}, T_{sb}^{[1]} \le u_{sb}^{[1]}, \ldots \nonumber\\
\ldots, \Big| T_{sum}^{[Q]} \Big| \le u_{sum}^{[Q]}, T_{lb}^{[Q]} \le u_{lb}^{[Q]}, T_{sb}^{[Q]} \le u_{sb}^{[Q]} \Big) - \mathbf{P} \Bigg( \Bigg| \sum_{k=1}^{N_{lb}^{[1]}} \zeta_{k}^{[1]}\Big({K_{sb}^{[1]}}+2\Big) \Bigg| \le u_{sum}^{[1]}, \nonumber \\
N_{lb}^{[1]} \sum_{k=1}^{N_{lb}^{[1]}} \Big(\zeta_k^{[1]}\Big(K_{sb}^{[1]}+1\Big)\Big)^2 \le u_{lb}^{[1]}, L_{sb}^{[1]} \sum_{j=0}^{K_{sb}^{[1]}} \Bigg(\sum_{k=1}^{N_{lb}^{[1]}} \zeta_k^{[1]}(j) \Bigg)^2 \le u_{sb}^{[1]}, \ldots \nonumber \\
\ldots, \Bigg| \sum_{k=1}^{N_{lb}^{[Q]}} \zeta_{k}^{[Q]}\Big(K_{sb}^{[Q]}+2\Big) \Bigg| \le u_{sum}^{[Q]}, N_{lb}^{[Q]} \sum_{k=1}^{N_{lb}^{[Q]}} \Big(\zeta_k^{[Q]}\Big(K_{sb}^{[Q]}+1\Big)\Big)^2 \le u_{lb}^{[Q]}, \nonumber \\
L_{lb}^{[Q]} \sum_{j=0}^{K_{sb}^{[Q]}} \Bigg(\sum_{k=1}^{N_{lb}^{[Q]}} \zeta_k^{[Q]}(j) \Bigg)^2 \le u_{sb}^{[Q]} \Bigg) \Bigg| \le \nonumber \\
\le \frac{\delta_{SH}^*\Big(N_{lcm},h,n, u_{sum}^{[1]}, u_{lb}^{[1]},u_{sb}^{[1]}, \ldots, u_{sum}^{[Q]}, u_{lb}^{[Q]},u_{sb}^{[Q]}\Big) + \delta_{BE}^*(N_{lcm},h)}{\sqrt{\underline{n}(N_{lcm},h,n)}} \le \nonumber \\
\le \frac{\overline{\delta}_{SH}^*(N_{lcm},h) + \delta_{BE}^*(N_{lcm},h)}{\sqrt{\underline{n}(N_{lcm},h,n)}}.\label{flaTh2}
\end{gather}
\end{theorem}
Let the conditions of Theorem \ref{th2Berry} be satisfied. Then the functions $f_{sum}^{[q]}$ and $f_{lb}^{[q]}$ are bounded for $1 \le q \le Q$ and $0 \le \delta_{BE}^*(N_{lcm},h)< \infty$ by the definition of $\delta_{BE}^*(N_{lcm},h)$. From the conditions of Theorem \ref{th2Berry} it follows that $N_{lcm} | N$ and therefore $N \ge N_{lcm}$. Substituting $N_{lcm}$ for $N$ in Lemma 2 \cite{Savelov13}, we obtain that if $N_{lb}^{[q]} \ge 2$ and $L_{sb}^{[q]} \ge 2$ for $1 \le q \le Q$, then the quantity $\overline{\delta}_{SH}^*(N_{lcm},h)$ is finite. Thus, the numerator of the fraction on the right-hand side of the formula \eqref{flaTh2} is finite in a wide class of cases. Since the quantity $\underline{n}(N_{lcm},h,n)=N_{lcm} h \cdot \big[\frac{n -N_{lcm}(h-1) }{N_{lcm}h}\big]$ satisfies the relation $\underline{n}(N_{lcm},h,n) \sim n$ (as $n \to \infty$), then the relation \eqref{flaTh2} is an analogue of the Berry-Esseen inequality.

If $m_{sum}^{[q]} = m_{lb}^{[q]}=1$, then the limit joint distribution of the statistics $T_{sum}^{[1]} , T_{lb}^{[1]}, T_{sb}^{[q]}, \ldots, T_{sum}^{[Q]} , T_{lb}^{[Q]}, T_{sb}^{[Q]}$ under $H_0$ can be found using Corollary \ref{thOsovmRaspr3stat} (setting the value $s$ from the condition of Corollary \ref{thOsovmRaspr3stat} equal to $h$ and substituting $J=0$) explicitly in terms of the quantities $\zeta_i^{[q]}, 1 \le i \le N_{lb}^{[q]}$; Theorem \ref{th2Berry} characterizes the rate of convergence to this limit distribution.

The proof of Theorem \ref{th2Berry} is obtained by almost verbatim repetition of the proof given in \cite{Savelov13}. We do not present the full proof because it is cumbersome.

\subsection{Asymptotic Independence Criterion for Four Types of Statistics} \label{zamechOTomChtoVseOshchaetsya}

As in case of the paper \cite{Savelov13}, the most important results of the paper \cite{Savelov12} (see Section ''Main Results'' in \cite{Savelov12}) also rely on \cite{Savelov11} and are also extended from the case of $\mathfrak{X} = \{0,1,\ldots,R-1\}$ to the case of an arbitrary $\mathfrak{X}$ with natural variations.

Let $T_{ij} \in \{ T_{sum}^{[1]}, T_{lb}^{[1]}, T_{sb}^{[1]}, \ldots, T_{sum}^{[Q]}, T_{lb}^{[Q]}, T_{sb}^{[Q]}, T_{quad}^{[1]},\ldots, T_{quad}^{[J]} \}$ for $1 \le i \le l$ and $1 \le j \le r_i$. Let us construct random variables $A_{T_{11}, T_{12}, \ldots, T_{1r_1}}, A_{T_{21}, T_{22}, \ldots, T_{2r_2}}, \ldots, A_{T_{l1}, T_{l2}, \ldots, T_{lr_l}}$ from random variables $\zeta_i^{[q]}$ in the same way as in \cite{Savelov12}. By asymptotic independence of statistics we mean the same notion as in \cite{Savelov11} and \cite{Savelov12}.

\begin{theorem} \label{thMainAsNez}
Let all the conditions of Theorem \ref{th2} be satisfied and the vector
$$\bigg(\eta_1^*(0), \ldots, \eta_1^*(K^*-1), \ldots, \eta_N^*(0), \ldots, \eta_N^*(K^*-1) \bigg)$$
be Gaussian. Let $T_{ij} \in \{ T_{sum}^{[1]}, T_{lb}^{[1]}, T_{sb}^{[1]}, \ldots, T_{sum}^{[Q]}, T_{lb}^{[Q]}, T_{sb}^{[Q]}, T_{quad}^{[1]},\ldots, T_{quad}^{[J]} \}$ for $1 \le i \le l$ and $1 \le j \le r_i$. The following conditions are equivalent:\\
a) the vectors $(T_{11}, T_{12}, \ldots, T_{1r_1}), (T_{21}, T_{22}, \ldots, T_{2r_2}), \ldots, (T_{l1}, T_{l2}, \ldots, T_{lr_l})$ are asymptotically independent;\\
b) the components of different vectors $(T_{11}, T_{12}, \ldots, T_{1r_1}), (T_{21}, T_{22}, \ldots, T_{2r_2}), \ldots$ $\ldots, (T_{l1}, T_{l2}, \ldots, T_{lr_l})$ are pairwise asymptotically independent, i.e. for any pair $i,j$ of different elements of the set $\{1,2,\ldots,l\}$, the statistics $T_{iu}$ and $T_{jv}$ are asymptotically independent for all $u \in \{1,2,\ldots,r_i\}, v \in \{1,2,\ldots,r_j\}$; \\
c) the vectors $A_{T_{11}, T_{12}, \ldots, T_{1r_1}}, A_{T_{21}, T_{22}, \ldots, T_{2r_2}}, \ldots, A_{T_{l1}, T_{l2}, \ldots, T_{lr_l}}$ are pairwise uncorrelated.
\end{theorem}

The proof of Theorem \ref{thMainAsNez} is similar to the proof of Theorem 3 \cite{Savelov13}.

Next, put $$G^* = \frac{1}{Nh}\Phi^*.$$
Note that $G^* = \mathbf{D}_{H_0}\eta_k^* = \frac{1}{Nh}\Phi^*$ for $1 \le k \le N$ (see Example \ref{zamechPriH0Horosho2}). We will number the elements of the matrix $G^*$ starting from zero. In other words,
\begin{equation} G_{i,j}^* = \mathrm{cov}_{H_0} (\eta_k^*(i), \eta_k^*(j)), \quad 0 \le i, j < K^*.\label{defGij}
\end{equation}
(cf. \cite[p.87]{Savelov12}). Note that the right-hand side of the equality \eqref{defGij} does not depend on \mbox{$k \in \{1, \ldots, N\}$}.

Let $T \in \{ T_{sum}^{[1]}, T_{lb}^{[1]}, T_{sb}^{[1]}, \ldots, T_{sum}^{[Q]}, T_{lb}^{[Q]}, T_{sb}^{[Q]} \}$ for $1 \le i \le l$ and $1 \le j \le r_i$. We define the set $\alpha_T$ in the same way as in \cite{Savelov12}.

\begin{corollary} \label{sledTeor3AsNez}
Assume that all conditions of Theorem \ref{th2} are satisfied. Let vectors $\eta_1^*, \ldots, \eta_N^*$ be independent, have a normal distribution (with, generally speaking, different expectations), and let the covariance matrix of each of them be equal to $G^*$. Let $T_{ij} \in \{ T_{sum}^{[1]}, T_{lb}^{[1]}, T_{sb}^{[1]}, \ldots, T_{sum}^{[Q]}, T_{lb}^{[Q]}, T_{sb}^{[Q]} \}$ for $1 \le i \le l$ and $1 \le j \le r_i$. The vectors $(T_{11}, T_{12}, \ldots, T_{1r_1}), (T_{21}, T_{22}, \ldots, T_{2r_2}), \ldots, (T_{l1}, T_{l2}, \ldots, T_{lr_l})$ are asymptotically independent if and only if for any pair $i,j$ of different elements of the set $\{1,2,\ldots,l\}$ the equality $G_{u,v}^*=0$ holds for all $u \in (\alpha_{T_{i1}} \cup \alpha_{T_{i2}} \cup \ldots \cup \alpha_{T_{ir_i}})$ and $v \in (\alpha_{T_{j1}} \cup \alpha_{T_{j2}} \cup \ldots \cup \alpha_{T_{jr_j}})$.
\end{corollary}
The proof of Corollary \ref{sledTeor3AsNez} is similar to the proof of Corollary 2 \cite{Savelov13}.

Moreover, Lemma 1 \cite{Savelov12} remains valid if in the formulation of its first statement the condition ''the conditions of Theorem 1 \cite{Savelov12} are satisfied'' is replaced by the condition ''all conditions of Theorem \ref{th2} (of the present paper) are satisfied and the vector
$$\bigg(\eta_1^*(0), \ldots, \eta_1^*(K^*-1), \ldots, \eta_N^*(0), \ldots, \eta_N^*(K^*-1) \bigg)$$ is Gaussian'', and in the formulation of its second statement the condition ''${\boldsymbol C}^{(N,s)} =(C(s,h), C(s,h), \ldots, C(s,h))$'' is replaced by the condition ''vectors $\eta_1^*, \ldots, \eta_N^*$ are independent and normally distributed, each with covariance matrix  $G^*$ (the expectations of the vectors $\eta_1^*, \ldots, \eta_N^*$ may be different)''. The proofs of the generalizations of the results of \cite{Savelov12} considered here are similar to the proofs of the corresponding statements in \cite{Savelov12}. We omit the full statements and proofs because they are cumbersome.

\begin{example} \label{exampleTh4}
Let us give an example of a situation where the conditions of Theorem \ref{thMainAsNez} are satisfied. Let $T_{ij} \in \{ T_{sum}^{[1]}, T_{lb}^{[1]}, T_{sb}^{[1]}, \ldots, T_{sum}^{[Q]}, T_{lb}^{[Q]}, T_{sb}^{[Q]}, T_{quad}^{[1]},\ldots, T_{quad}^{[J]} \}$ for $1 \le i \le l$ and $1 \le j \le r_i$. Let the parameters of the statistics $T_{quad}^{[1]},\ldots, T_{quad}^{[J]}$ be fixed and the conditions of Corollary 3 \cite{Savelov11} are satisfied. Then all conditions of Theorem 2 \cite{Savelov11} are satisfied, and by virtue of the considerations from Example \ref{zamechPriH0Horosho3} the relation
\eqref{XStarshoditsyaKetaStar} is satisfied for the corresponding vectors $\eta_1^*, \ldots, \eta_N^*$ from Theorem 2 \cite{Savelov11}, which, by virtue of Corollary 3 \cite{Savelov11}, are independent vectors, each of which has a normal distribution with covariance matrix $G^*$. The function $f_{sum}^{[q]}$ is bounded, since it is defined on a finite set $\mathfrak{X}^{m_{sum}^{[q]}}$, where $\mathfrak{X} = \{0,1,\ldots,R-1\}$. Similarly, we get that the function $f_{lb}^{[q]}$ is also bounded. This means that the conditions of Theorem \ref{thMainAsNez} are satisfied. We also note that examples of Markov alternatives $H_1$ for which the conditions of Corollary 3 \cite{Savelov11} are satisfied are considered in \cite{Savelov11} (see hypotheses $\tilde{H}_1$ and $\hat{H}_1$).
\end{example}

\begin{example}
Let us give an example of a situation where the conditions of Corollary \ref{sledTeor3AsNez} are satisfied. Let $T_{ij} \in \{ T_{sum}^{[1]}, T_{lb}^{[1]}, T_{sb}^{[1]}, \ldots, T_{sum}^{[Q]}, T_{lb}^{[Q]}, T_{sb}^{[Q]} \}$ for $1 \le i \le l$ and $1 \le j \le r_i$. Suppose that the parameters of the statistics $T_{quad}^{[1]},\ldots, T_{quad}^{[J]}$ are fixed and the conditions of Corollary 3 \cite{Savelov11} are satisfied. Reasoning in the same way as in Example \ref{exampleTh4}, we find that the conditions of \mbox{Corollary \ref{sledTeor3AsNez}} are satisfied.
\end{example}

\section{Proofs}

\subsection{Proof of Theorem \ref{th1}}

First, we prove auxiliary results without assuming that relation \eqref{XshoditsyaKeta} holds.

From the boundedness of $f_{sum}$ and $f_{lb}$ it follows that if $j_n \in \{0,\ldots,n-m_{sum}\}$ for $n \ge 1$, then
\begin{equation}
\frac{f_{sum}(\varepsilon_{j_n+1}, \ldots, \varepsilon_{j_n+m_{sum}})}{\sqrt{n}} = o_P(1), \quad n \to \infty, \label{prenebrejimMalostSlag1}
\end{equation}
and if $j_n \in \{0,\ldots,n-m_{lb}\}$ for $n \ge 1$, then
\begin{equation}
\frac{f_{lb}(\varepsilon_{j_n+1}, \ldots, \varepsilon_{j_n+m_{lb}})}{\sqrt{n}} = o_P(1), \quad n \to \infty. \label{prenebrejimMalostSlag2}
\end{equation}
Taking into account \eqref{prenebrejimMalostSlag2}, for $k \in \{1,\ldots, N_{lb}\}$ we obtain that
\begin{gather}
\mathcal{Y}_k^{(N,s,h)}({K_{sb}}+1)  =  \sum_{t = \frac{N}{N_{lb}}(k-1) + 1}^{\frac{N}{N_{lb}}k} \mathcal{X}_t^{(N,s,h)}({K_{sb}}+1) = \nonumber \\
=\! \! \! \! \! \! \!\sum_{t = \frac{N}{N_{lb}}(k-1) + 1}^{\frac{N}{N_{lb}}k} \! \! \! \! \! \! n^{-\frac12} \! \! \!  \! \! \! \sum_{\frac{L(t-1)}h\le i \le \frac{Lt-s}h} \!\!\!\!\! \big( f_{K_{sb}+1}^{(s)}(\varepsilon_{hi+1}, \ldots, \varepsilon_{hi+s})-
\mathbf{E}_{H_0} f_{K_{sb}+1}^{(s)}(\varepsilon_{1}, \ldots, \varepsilon_{s})\big) \!= \nonumber\\
 = \! \! \! \! \! \! \! \sum_{t = \frac{N}{N_{lb}}(k-1) + 1}^{\frac{N}{N_{lb}}k}   \sum_{\frac{L(t-1)}h\le i \le \frac{Lt-s}h} \frac{1}{\sqrt{n}}
 \sum_{u=1}^{h} \frac{ f_{lb}(\varepsilon_{hi+u}, \ldots, \varepsilon_{hi+u+{m_{lb}}-1})\!-\!
\mathbf{E}_{H_0} f_{lb}(\varepsilon_1,  \ldots, \varepsilon_{m_{lb}})}{\sigma_{lb}}
 = \nonumber \\
 = \! \! \! \! \! \! \! \sum_{t = \frac{N}{N_{lb}}(k-1) + 1}^{\frac{N}{N_{lb}}k} \sum_{L(t-1)+1\le j \le Lt-s+h}   \! \! \! \! \! \! \! \! \!   \! \! \! \! \! \! \! \frac{
f_{lb}(\varepsilon_{j},\!\ldots,\!\varepsilon_{j+{m_{lb}}-1})  - \mathbf{E}_{H_0}  f_{lb}(\varepsilon_1,\!\ldots,\!\varepsilon_{{m_{lb}}})
  }{\sqrt{n}\sigma_{lb}} + o_P(1)
 = \nonumber \\
  = \!\!\! \sum_{j=\frac{LN}{N_{lb}} (k-1)+1}^{ \frac{LN}{N_{lb}} k-{m_{lb}}+1}  \!\!\! \frac{
f_{lb}(\varepsilon_{j},  \ldots, \varepsilon_{j+{m_{lb}}-1})  - \mathbf{E}_{H_0}  f_{lb}(\varepsilon_1,  \ldots, \varepsilon_{{m_{lb}}})
  }{\sqrt{n}\sigma_{lb}} + o_P(1)
 = \nonumber \\
=\frac{  \sum_{j={L_{lb}}(k-1)+1}^{{L_{lb}}k-{m_{lb}}+1} \big(f_{lb}(\varepsilon_{j}, \ldots, \varepsilon_{j+{m_{lb}}-1})  - \mathbf{E}_{H_0}  f_{lb}(\varepsilon_1, \ldots, \varepsilon_{{m_{lb}}})  \big)}{\sqrt{n}\sigma_{lb}}  + o_P(1) = \nonumber \\
= \frac{W_k - (L_{lb}-m_{lb}+1)E_{lb}}{\sqrt{n}\sigma_{lb}} + o_P(1) \nonumber
\end{gather}
and therefore
\begin{gather}
\frac{W_k - (L_{lb}-m_{lb}+1)E_{lb}}{\sqrt{n}\sigma_{lb}}  = \mathcal{Y}_k^{(N,s,h)}({K_{sb}}+1) + o_P(1).\label{equivDlyaWk}
\end{gather}
By similar considerations, using the relation \eqref{prenebrejimMalostSlag1} we find that
$$\sum_{k=1}^{N} \mathcal{X}_k^{(N,s,h)}({K_{sb}}+2) = \sum_{i=1}^{n-m_{sum}+1}\frac{f_{sum}(\varepsilon_i, \ldots, \varepsilon_{i+m_{sum}-1})\,-E_{sum}}{\sqrt{n} \sigma_{sum}} +o_p(1)$$
and, as a corollary,
\begin{equation} \label{flaTsumAsRazlvariant1}
\begin{aligned}
\frac{\sqrt{n-m_{sum}+1}}{\sqrt{n}}T_{sum}\!=\!\!\sum_{i=1}^{n-m_{sum}+1}\frac{f_{sum}(\varepsilon_i, \ldots, \varepsilon_{i+m_{sum}-1})\,-E_{sum}}{\sqrt{n} \sigma_{sum}}\,= \\
=  \sum_{k=1}^{N} \mathcal{X}_k^{(N,s,h)}({K_{sb}}+2) +o_p(1) = \sum_{k=1}^{N_{lb}} \mathcal{Y}_k^{(N,s,h)}({K_{sb}}+2) +  o_p(1).
\end{aligned}
\end{equation}
For $j \in \{0,1,\ldots, K_{sb}\}$
\begin{gather}
 \sum_{k=1}^{N_{lb}} \mathcal{Y}_k^{(N,s,h)}(j) =  \sum_{k=1}^{N_{lb}} \sum_{i = 1}^{\frac{N}{N_{lb}}} \mathcal{X}_{\frac{N}{N_{lb}}(k-1)+i}^{(N,s,h)}(j) = \sum_{k = 1}^N \mathcal{X}_k^{(N,s,h)}(j) = \nonumber \\
=\!\sum_{k = 1}^N
  n^{-\frac12}\!\!\!\!\!\!\!\!\!\sum_{\frac{L(k-1)}h\le i \le \frac{Lk-s}h} \!\!\! \Big(  f_j^{(s)}(\varepsilon_{hi+1}, \ldots, \varepsilon_{hi+s}) -  \mathbf{E}_{H_0}f_j^{(s)}(\varepsilon_{1},  \ldots, \varepsilon_{s})  \Big) =
\nonumber \\
=\! \sum_{k = 1}^N
 \sum_{\frac{L(k-1)}h\le i \le \frac{Lk-s}h}\!
 \sum_{u=1}^{h / L_{sb}} \! \frac{I_{f_{sb}(\varepsilon_{hi+{L_{sb}}(u-1)+1}, \ldots, \varepsilon_{hi+{L_{sb}}u})\in \alpha_{sb}(j) } - E_{sb}(j)}{\sqrt{n E_{sb}(j)}} \!=
\nonumber\\
= \! \sum_{k = 1}^N
 \sum_{\frac{L}{L_{sb}}(k-1)+1 \le t \le \frac{Lk}{L_{sb}}} \!\!\!\!\!\!\!\!\!
 \frac{I_{f_{sb}(\varepsilon_{L_{sb}(t-1)+1}, \ldots, \varepsilon_{L_{sb}t}) \in \alpha_{sb}(j) }-E_{sb}(j)}{ \sqrt{ n E_{sb}(j)}} +o_p(1) =
\nonumber\\
=  \sum_{1 \le t \le \frac{LN}{L_{sb}}}
 \frac{I_{f_{sb}(\varepsilon_{L_{sb}(t-1)+1}, \ldots, \varepsilon_{L_{sb}t}) \in \alpha_{sb}(j) }-E_{sb}(j)}{ \sqrt{ n E_{sb}(j)}} +o_p(1) =
\nonumber\\
= \frac{\sum_{t=1}^{{N_{sb}}} \big(I_{f_{sb}(\varepsilon_{L_{sb}(t-1)+1}, \ldots, \varepsilon_{L_{sb}t}) \in \alpha_{sb}(j) }-E_{sb}(j)\big)}{ \sqrt{ n E_{sb}(j)}} +o_p(1)
= \frac{ w(j) - {N_{sb}}E_{sb}(j)}{ \sqrt{ n E_{sb}(j)}} +o_p(1). \nonumber
\end{gather}
From here and from \eqref{equivDlyaWk} we find that
\begin{gather}
 \frac{1}{\sqrt{n}} \cdot\bigg( \frac{W_1 - (L_{lb}-m_{lb}+1)E_{lb}}{\sigma_{lb}},
\ldots, \frac{W_{N_{lb}} - (L_{lb}-m_{lb}+1)E_{lb}}{\sigma_{lb}},  \nonumber \\
\frac{ w(0) - {N_{sb}}E_{sb}(0)}{\sqrt{E_{sb}(0)}}, \ldots,  \frac{w(K_{sb}) - {N_{sb}} E_{sb}(K_{sb})  }{\sqrt{ E_{sb}(K_{sb})  }} \bigg) =  \nonumber \\
=  \bigg( \mathcal{Y}_1^{(N,s,h)}({K_{sb}}+1), \ldots, \mathcal{Y}_{N_{lb}}^{(N,s,h)}({K_{sb}}+1), \nonumber \\
\sum_{k=1}^{N_{lb}} \mathcal{Y}_k^{(N,s,h)}(0), \ldots, \sum_{k=1}^{N_{lb}} \mathcal{Y}_k^{(N,s,h)}({K_{sb}})\bigg) +o_p(1).\label{flaTlbsbAsRazl}
\end{gather}
Let us show that \eqref{flaTsumAsRazlvariant1} and the boundedness of $f_{sum}$ imply the relation \eqref{flaTsumAsRazl}. Indeed, if $|f_{sum}| \le C$ for some $C \in
\mathbb{R}$, then $|E_{sum}| \le C$ and $$\Bigg|\sum_{i=1}^{n-m_{sum}+1} \frac{f_{sum}(\varepsilon_i, \ldots, \varepsilon_{i+m_{sum}-1}) - E_{sum}}{n}\Bigg| \le \frac{2C(n-m_{sum}+1)}{n} \le 2C.$$
From this and from \eqref{flaTsumAsRazlvariant1} we obtain that the value $\frac{\sqrt{n-m_{sum}+1}}{n}T_{sum}$ (and, as a consequence, the value $\frac{T_{sum}}{\sqrt{n}}$) is bounded. From \eqref{flaTsumAsRazlvariant1} and the boundedness of the value $\frac{T_{sum}}{\sqrt{n}}$ it follows that
\begin{gather}
T_{sum} =\sqrt{1 - \frac{m_{sum}-1}n} T_{sum} + \bigg(1 - \sqrt{1 - \frac{m_{sum}-1}n}\bigg)T_{sum} = \nonumber \\
= \sqrt{1 - \frac{m_{sum}-1}n}T_{sum} + o_P(1) = \sum_{k=1}^{N_{lb}} \mathcal{Y}_k^{(N,s,h)}({K_{sb}}+2) +o_p(1). \nonumber
\end{gather}
The equality \eqref{flaTsumAsRazl} is proved.

Further, assume that the relation \eqref{XshoditsyaKeta} holds. From \eqref{defYk} and \eqref{defZetaK} it follows that
\begin{gather}\Big(\mathcal{Y}_1^{(N,s,h)}, \mathcal{Y}_2^{(N,s,h)}, \ldots, \mathcal{Y}_{N_{lb}}^{(N,s,h)}\Big)
 =\sum_{i = 1}^{\frac{N}{N_{lb}}} \Big( \mathcal{X}_{i}^{(N,s,h)}, \mathcal{X}_{\frac{N}{N_{lb}}+i}^{(N,s,h)}, \ldots, \mathcal{X}_{\frac{N}{N_{lb}}(N_{lb}-1)+i}^{(N,s,h)}
 \Big), \label{Y4erezX} \\
(\zeta_1, \zeta_2, \ldots, \zeta_{N_{lb}}) = \sum_{i = 1}^{\frac{N}{N_{lb}}}\Big( \eta_{i}, \eta_{\frac{N}{N_{lb}}+i}, \ldots, \eta_{\frac{N}{N_{lb}}(N_{lb}-1)+i}\Big). \label{zeta4erezEta}
\end{gather}
Taking into account \eqref{XshoditsyaKeta}, \eqref{Y4erezX}, and \eqref{zeta4erezEta}, we obtain the relation \eqref{limY}.

Further, from \eqref{limY} it follows that the sequence of random vectors
\begin{gather}\bigg(\mathcal{Y}_1^{(N,s,h)}({K_{sb}}+1), \mathcal{Y}_2^{(N,s,h)}({K_{sb}}+1), \ldots, \mathcal{Y}_{N_{sb}}^{(N,s,h)}({K_{sb}}+1), \sum_{i=1}^{N_{sb}} \mathcal{Y}_i^{(N,s,h)}(0), \ldots, \sum_{i=1}^{N_{sb}} \mathcal{Y}_i^{(N,s,h)}({K_{sb}})\bigg)\nonumber
\end{gather}
converges in distribution as $n \to \infty$, from which, by \eqref{flaTlbsbAsRazl}, we obtain that
 \begin{gather}
\!\bigg( \frac{L_{lb}}{n} T_{lb}, \frac{{N_{sb}}}{n} T_{sb}  \bigg) \!\! = \!\! \Bigg( \sum_{k=1}^{N_{lb}} \big(\mathcal{Y}_k^{(N,s,h)}(K_{sb}+1)\big)^2\!, \sum_{j=0}^{K_{sb}} \bigg(\sum_{k=1}^{N_{lb}} \mathcal{Y}_k^{(N,s,h)}(j) \bigg)^2 \Bigg)\! +\! o_P(1). \nonumber
\end{gather}
This means that the equalities \eqref{flaTlbAsRazl} and \eqref{flaTsbAsRazl} are true.

It follows from \eqref{flaTsumAsRazl}, \eqref{flaTlbAsRazl}, \eqref{flaTsbAsRazl} and \eqref{limY} that \eqref{th1LimDistr} holds. Theorem \ref{th1} is proved.

\begin{remark} \label{MojnoNeogr}
Note that the boundedness condition for the functions $f_{sum}$ and $f_{lb}$ from Theorem \ref{th1} is only necessary to ensure that the equalities \eqref{prenebrejimMalostSlag1}--\eqref{prenebrejimMalostSlag2} hold and that the equality \eqref{flaTsumAsRazl} can be derived from the relation \eqref{flaTsumAsRazlvariant1}. Let's show that if $H_0$ holds, then there is no need to assume the boundedness of the functions $f_{sum}$ and $f_{lb}$. Indeed, if $H_0$ is true, then $\frac{f_{sum}(\varepsilon_{j_n+1}, \ldots, \varepsilon_{j_n+m_{sum}})}{\sqrt{n}}$ coincides in distribution with the random variable $\frac{f_{sum}(\varepsilon_1, \ldots, \varepsilon_{m_{sum}})}{\sqrt{n}}$, which tends to zero in probability \Big (and similarly $\frac{f_{lb}(\varepsilon_{j_n+1}, \ldots, \varepsilon_{j_n+m_{lb}})}{\sqrt{n}} = o_P(1)$\Big), and the relation \eqref{flaTsumAsRazl} follows from the relation \eqref{flaTsumAsRazlvariant1}, since the expansion $T_{sum} =\sqrt{\frac{n-m_{sum}+1}n} T_{sum} + \bigg(1 - \sqrt{1 - \frac{m_{sum}-1}n}\bigg)T_{sum}$ holds, in which the second term is the product of an infinitesimal sequence and a convergent sequence (see \eqref{NotJointDistrOfEvery}). Note also that in the case where the elements of the tested sequence take values in a fixed finite set (this situation was considered in \cite{Savelov11, Savelov12, Savelov13}), the boundedness condition for the functions $f_{sum}$ and $f_{lb}$ is satisfied automatically.
\end{remark}

\subsection{Proof of Corollary \ref{corolOsovmRaspr3stat}}

If the functions $f_{sum}$ and $f_{lb}$ are bounded, then the statement of Corollary \ref{corolOsovmRaspr3stat} is a direct consequence of Theorem \ref{th1} and Example \ref{zamechPriH0Horosho1}. The case when $f_{sum}$ or $f_{lb}$ is unbounded is reduced to the previous one by Remark \ref{MojnoNeogr}.

\subsection{Proof of Theorem \ref{th2}}
Let $1 \le q \le Q$. We suppose that $s \ge \underline{s}^*$. Therefore, $s \ge h + \max\Big(m_{sum}^{[q]},m_{lb}^{[q]}\Big) -1$. Applying Theorem \ref{th1} to the triple of statistics $(T_{sum}, T_{lb}, T_{sb})=\Big(T_{sum}^{[q]}, T_{lb}^{[q]} , T_{sb}^{[q]}\Big)$, we obtain the equalities \eqref{flaTsumAsRazl2}, \eqref{flaTlbAsRazl2}, \eqref{flaTsbAsRazl2}.

From the definition of the vectors $\zeta_k^{[q]}$ $\big(1 \le k \le N_{lb}^{[q]}, \ 1 \le q \le Q \big)$ it follows that there exists a linear mapping $A$ such that
$$\Big(\zeta_1^{[1]}, \ldots, \zeta_{N_{lb}^{[1]}}^{[1]}, \ldots, \zeta_1^{[Q]}, \ldots, \zeta_{N_{lb}^{[Q]}}^{[Q]}\Big) = A (\eta_1^*, \ldots, \eta_N^*).$$
A straightforward check shows that \begin{gather}
\Big(\mathcal{Y}_1^{[1],(N,s,h)}, \ldots, \mathcal{Y}_{N_{lb}^{[1]}}^{[1],(N,s,h)}, \ldots, \mathcal{Y}_1^{[Q],(N,s,h)}, \ldots, \mathcal{Y}_{N_{lb}^{[Q]}}^{[Q],(N,s,h)}\Big) =A\Big(\mathcal{X}_1^{*(N,s,h)}, \ldots, \mathcal{X}_N^{*(N,s,h)}\Big). \nonumber
\end{gather}
This means that \eqref{limYbig} follows from \eqref{XStarshoditsyaKetaStar}.

By \eqref{flaTsumAsRazl2} and \eqref{limYbig}, for $1 \le q \le Q$, the equality $T_{sum}^{[q]} = O_P(1)$ holds as $n \to \infty$ (i.e., $T_{sum}^{[q]}$ is bounded in probability). This, together with \eqref{deTquadJ} and \eqref{flaTsumAsRazl2}, implies the equality \eqref{flaTquadAsRazl}.

It follows from \eqref{flaTsumAsRazl2}, \eqref{flaTlbAsRazl2}, \eqref{flaTsbAsRazl2}, \eqref{flaTquadAsRazl}, and \eqref{limYbig} that \eqref{th2LimDistrGeneral} holds. Theorem \ref{th2} is proved.

Note that Theorem \ref{th2} also holds in the case where quadratic statistics are absent (i.e., $J=0$), if for $J=0$ in the formulation of Theorem \ref{th2} the vector \eqref{th2LimDistrGeneral} lacks components of the form $\theta^{[j]}$, and the equality \eqref{flaTquadAsRazl} is also absent. Indeed, for the case $J=0$, the equalities \eqref{flaTsumAsRazl2}, \eqref{flaTlbAsRazl2}, \eqref{flaTsbAsRazl2}, and the relation \eqref{limYbig} have already been proven, and from them an analogue of the relation \eqref{th2LimDistrGeneral} follows, in which components of the form $\theta^{[j]}$ (and quadratic statistics) are absent.

\begin{remark} \label{MojnoNeogr2}
As the analysis of the proof of Theorem \ref{th2} shows, if $H_0$ is true, then by virtue of Remark \ref{MojnoNeogr} the statement of Theorem \ref{th2} will remain valid even if the boundedness condition of the functions $f_{sum}^{[q]}$ and $f_{lb}^{[q]}$ ($1 \le q \le Q$) is not satisfied (but the remaining conditions of Theorem \ref{th2} are satisfied).
\end{remark}

\subsection{Proof of Corollary \ref{thOsovmRaspr3stat}}

Applying Theorem \ref{th2} and taking into account Example \ref{zamechPriH0Horosho2}, we obtain the statement of Corollary \ref{thOsovmRaspr3stat} in the case where we additionally assume that the functions $f_{sum}^{[q]}$ and $f_{lb}^{[q]}$ are bounded for $1 \le q \le Q$. The case where not all functions $f_{sum}^{[q]}$ and $f_{lb}^{[q]}$ are bounded is reduced to the previous one via Remark \ref{MojnoNeogr2}.

\subsection{Proof of Lemma \ref{LemmaProNeBlizkH1}}

The proof of Lemma \ref{LemmaProNeBlizkH1} is similar to the proof of Lemma 5 \cite{Savelov11}.

\begin{remark} \label{ProLemmyProNeblizkH1}
As an analysis of the proof of Lemma \ref{LemmaProNeBlizkH1} shows (cf. Remark 6 \cite{Savelov11}), if convergence in probability is replaced by almost-sure convergence throughout its formulation, then the statement of Lemma \ref{LemmaProNeBlizkH1} remains true.
\end{remark}

\subsection{On the statistic $T_{SO}$ from Example \ref{primersmarsaSerialOver}}\label{oStatistTSo}

Put $$R = 2^{r_{bits}}.$$ From the sequence $\varepsilon_1, \ldots, \varepsilon_n$, we construct a sequence of binary vectors $\varkappa_1, \ldots, \varkappa_n$, where $\varkappa_i = g_{bits}(\varepsilon_i)$. The random vector $\varkappa_i$ takes values from the set $\{k_0, \ldots, k_{R-1}\}$, where $k_0 = \underbrace{(0,0,\ldots,0)}_{\text{$r_{bits}$ times}}, \ldots, k_{R-1} = \underbrace{(1,1,\ldots,1)}_{\text{$r_{bits}$ times}}$. Let $m \ge 2$ and $(i_1,\ldots, i_m) \in \{k_0, k_1, \ldots, k_{R-1}\}^m$. The number of occurrences of the chain $(i_1,\ldots, i_m)$ in the sequence $\varkappa_1, \varkappa_2, \ldots, \varkappa_{n-1}, \varkappa_n, \varkappa_{1}, \varkappa_{2}, \ldots, \varkappa_{m-1}$ is denoted by $\nu_{i_1,\ldots, i_m}$. In other words, $$\nu_{i_1,\ldots, i_m} = \sum_{i=1}^n I_{(\tilde{\varkappa}_i,\ldots, \tilde{\varkappa}_{i+m-1}) = (i_1,\ldots, i_m)},$$
where $\tilde{\varkappa}_i = \varkappa_i I_{i \le n} + \varkappa_{i-n} I_{i > n}$. The quantity $\nu_{i_1,\ldots, i_{m-1}}$ is defined similarly. By a sum of the form $\sum_{i_1,\ldots, i_d}$ we will henceforth mean the sum $\sum_{(i_1,\ldots, i_d) \in \{ k_0,k_1, \ldots, k_{R-1}\}^d}$. Put
$$ T_{SO} = \frac{R^{m}}n\sum_{i_1,\ldots,i_{m}} \Big( \nu_{i_1,\ldots, i_{m}} - \frac{n}{R^{m}}\Big)^2 - \frac{R^{m-1}}n\sum_{i_1,\ldots,i_{m-1}} \Big( \nu_{i_1,\ldots, i_{m-1}} - \frac{n}{R^{m-1}}\Big)^2.$$
This statistic is used in the ''smarsa\_SerialOver'' test of the TestU01 \cite{TestU01ShortGuide} package. The package \cite{TestU01ShortGuide} does not give an explicit formula, but it's given in \cite{Altman1988}, which \cite{TestU01ShortGuide} relies on. If $H_0^{U[0,1]}$ holds, then the random variables $\varkappa_i$ are independent and have a uniform distribution on the set $\{k_0, \ldots, k_{R-1}\}$, therefore (see \cite{Good1, Good2, Altman1988, Savelov8}) \begin{equation}
T_{SO} \overset{d}{\to} \chi_{R^m-R^{m-1}}^2. \label{limTSo}
\end{equation}
\begin{lemma} \label{lemmaProTSo}
Let $H_0^{U[0,1]}$ be true and the numbers $m \ge 2$ and $r_{bits}$ be fixed (independent of $n$). Then $T_{SO}$ coincides up to $o_P(1)$ with some quadratic statistic.
\end{lemma}
Let us prove Lemma \ref{lemmaProTSo}. For $(i_1, \ldots, i_m) \in \{ k_0,k_1,\ldots, k_{R-1} \}^m$ put
\begin{gather}
\mathcal{Z}_{i_1,\ldots, i_m} = \sqrt{n}\bigg( \frac{\nu_{i_1,\ldots, i_m}}{n} - R^{-m} \bigg), \quad \nu'_{i_1,\ldots, i_m} = \sum_{i=1}^{n-m+1} I_{(\varkappa_i,\ldots, \varkappa_{i+m-1}) = (i_1,\ldots, i_m)}, \nonumber \\ \mathcal{Z}'_{i_1,\ldots, i_m} = \sqrt{n-m+1}\bigg( \frac{\nu'_{i_1,\ldots, i_m}}{n-m+1} - R^{-m} \bigg).\nonumber
\end{gather}
Note that $\sum_{i_m} \nu_{i_1,\ldots, i_m} = \nu_{i_1,\ldots, i_{m-1}}$. Moreover, for any real numbers $c_0, \ldots, c_{R-1}$, the inequality $R\sum_{i=0}^{R-1}c_i^2 \ge \big(\sum_{i=0}^{R-1} c_i\big)^2$ holds. Hence,
\begin{gather}
T_{SO} = \frac{R^{m}}n\sum_{i_1,\ldots,i_{m-1}} \sum_{i_{m}}\Big( \nu_{i_1,\ldots, i_{m}} - \frac{n}{R^{m}}\Big)^2 -  \frac{R^{m-1}}n\sum_{i_1,\ldots,i_{m-1}} \bigg( \sum_{i_{m}}\Big( \nu_{i_1,\ldots, i_{m}} - \frac{n}{R^{m}} \Big) \bigg)^2 \ge 0, \nonumber
\end{gather}
i.e. $ T_{SO}$ is a non-negative definite quadratic form (with coefficients independent of $n$) of the random variables $\frac{1}{\sqrt{n}}(\nu_{i_1,\ldots, i_{m}} - \frac{n}{R^{m}}) = \mathcal{Z}_{i_1,\ldots, i_m}$. It is easy to write out explicitly the constants
$c_{i_1,\ldots,i_m,j_1,\ldots,j_m}$ such that $T_{SO} = \sum_{i_1,\ldots, i_m, j_1,\ldots, j_m} c_{i_1,\ldots,i_m,j_1,\ldots,j_m} \mathcal{Z}_{i_1,\ldots, i_m} \mathcal{Z}_{j_1,\ldots, j_m}$. Let $\tilde{T}_{SO}$ denote the same quadratic form in random variables $\mathcal{Z}'_{i_1,\ldots, i_m}$:
\begin{equation}\tilde{T}_{SO} = \sum_{i_1,\ldots, i_m, j_1,\ldots, j_m} c_{i_1,\ldots,i_m,j_1,\ldots,j_m} \mathcal{Z}'_{i_1,\ldots, i_m} \mathcal{Z}'_{j_1,\ldots, j_m}.\label{defTildeTSo}
\end{equation}
From the inequalities $|\nu_{i_1,\ldots, i_m} - \nu'_{i_1,\ldots, i_m}| \le m-1$ and $|\nu'_{i_1,\ldots, i_m}| \le n-m+1$ it follows that
\begin{equation}
\mathcal{Z}'_{i_1,\ldots, i_m} = \mathcal{Z}_{i_1,\ldots, i_m} + o_P(1), \quad n \to \infty.\label{ZshtrihAndz}
\end{equation}
Since $H_0^{U[0,1]}$ is true, the random variables $\varkappa_i$ are independent and have a uniform distribution on the set $\{k_0, \ldots, k_{R-1}\}$. By virtue of the central limit theorem for $(m-1)$-dependent random vectors, the limit distribution of an $R^m$-dimensional vector with components $\mathcal{Z}'_{i_1,\ldots, i_m}$, $(i_1, \ldots, i_m) \in \{k_0, k_1, \ldots, k_{R-1}\}^m$, is a multivariate normal distribution with zero mean. From this, from \eqref{ZshtrihAndz} and from the definition of the statistic $\tilde{T}_{SO}$ it follows that
\begin{equation}
\tilde{T}_{SO} = T_{SO} + o_P(1).\label{equivTSoAndTildeTSo}
\end{equation}
Taking into account \eqref{limTSo}, we obtain that
\begin{equation}
\tilde{T}_{SO} \overset{d}{\to} \chi_{R^m-R^{m-1}}^2. \label{limTildeTSo}
\end{equation}

Next, denote the variance of the limit distribution of the random variable $\mathcal{Z}'_{i_1,\ldots, i_m}$ by $\sigma_{i_1,\ldots,i_m}^2$. We assume that $\sigma_{i_1,\ldots,i_m} \ge 0$. The following convergence holds:
\begin{gather}
\mathcal{Z}'_{i_1,\ldots, i_m} \overset{d} {\to} \mathcal{N}(0,\sigma_{i_1,\ldots,i_m}^2), \quad n \to \infty. \label{shstZshtrih}
\end{gather}
Let us show that if $\sigma_{i_1,\ldots,i_m} \ne 0$, then $\frac{\mathcal{Z}'_{i_1,\ldots, i_m}}{\sigma_{i_1,\ldots,i_m}}$ is a summing statistic. Indeed, the summing statistic $T_{sum}$ (see \eqref{defTsum}) is constructed from the function $f_{sum}$ as follows. First, a non-negative value $\sigma_{sum}$ is calculated, the square of which is the variance of the limiting (as $n \to \infty$) distribution from the central limit theorem for $({m_{sum}}-1)$-dependent random variables $f_{sum}(\varepsilon_i, \ldots, \varepsilon_{i+{m_{sum}}-1})$, $i = 1, \ldots, n-m_{sum}+1$, i.e. $\sigma_{sum}^2$ is the limit variance of the sum $$\frac{ \sum_{i=1}^{n-m_{sum}+1} \big( f_{sum}(\varepsilon_i, \ldots, \varepsilon_{i+{m_{sum}}-1}) - \mathbf{E}_{H_0} f_{sum}(\varepsilon_1, \ldots, \varepsilon_{{m_{sum}}}) \big)}{\sqrt{n - m_{sum} +1}}.$$ In the case where $\sigma_{sum} \in (0,\infty)$, the statistic $T_{sum}$ is by definition set equal to the normalized sum $\frac{ \sum_{i=1}^{n-m_{sum}+1} \big( f_{sum}(\varepsilon_i, \ldots, \varepsilon_{i+{m_{sum}}-1}) - \mathbf{E}_{H_0} f_{sum}(\varepsilon_1, \ldots, \varepsilon_{{m_{sum}}}) \big)}{\sigma_{sum}\sqrt{n - m_{sum} +1}}$. Put $m_{sum} = m$ and $f_{sum}(\theta_1,\ldots,\theta_{m_{sum}}) = I_{g_{bits}(\theta_1) = i_1,\ldots, g_{bits}(\theta_{m_{sum}}) = i_{m_{sum}}}$.
Note that
\begin{gather}
\mathcal{Z}'_{i_1,\ldots, i_m} = \frac{ \sum_{i=1}^{n-m_{sum}+1} \big( f_{sum}(\varepsilon_i, \ldots, \varepsilon_{i+{m_{sum}}-1}) - \mathbf{E}_{H_0} f_{sum}(\varepsilon_1, \ldots, \varepsilon_{{m_{sum}}}) \big)}{\sqrt{n - m_{sum} +1}}.\nonumber
\end{gather}
Therefore, $\sigma_{sum}^2$ is the limit variance of $\mathcal{Z}'_{i_1,\ldots, i_m}$. By virtue of \eqref{shstZshtrih} $\sigma_{sum}^2 = \sigma_{i_1,\ldots,i_m}^2$, and the numbers $\sigma_{sum}$ and $\sigma_{i_1,\ldots,i_m}$ are non-negative by definition. From this and from the relation $\sigma_{i_1,\ldots,i_m} \ne 0$ it follows that $\sigma_{sum} = \sigma_{i_1,\ldots,i_m}$ and the quantity $\frac{\mathcal{Z}'_{i_1,\ldots, i_m}}{\sigma_{i_1,\ldots,i_m}} = \frac{\mathcal{Z}'_{i_1,\ldots, i_m}}{\sigma_{sum}}$ is a summing statistic.

If $\sigma_{i_1,\ldots,i_m} = 0$, then it follows from \eqref{shstZshtrih} that the equality $\mathcal{Z}'_{i_1,\ldots, i_m} = o_P(1)$ holds. Combining this with \eqref{defTildeTSo} and \eqref{limTildeTSo} we obtain that not all numbers $\sigma_{i_1,\ldots,i_m}$, $(i_1,\ldots, i_m) \in \{k_0, k_1, \ldots, k_{R-1}\}^m$, are equal to zero. Moreover, due to \eqref{defTildeTSo}
\begin{equation}\tilde{T}_{SO} = \!\!\!\!\sum_{i_1,\ldots, i_m, j_1,\ldots, j_m: \sigma_{i_1,\ldots,i_m} \sigma_{j_1,\ldots,j_m} \ne 0} \!\!\!\!\!\! \!\!\! \!\!\!\!\!\! \!\!\! c_{i_1,\ldots,i_m,j_1,\ldots,j_m} \mathcal{Z}'_{i_1,\ldots, i_m} \mathcal{Z}'_{j_1,\ldots, j_m} + o_P(1).\label{predstavlTSo}
\end{equation}
If some quadratic form $\sum_{i,j = 1}^d c_{i,j}x_ix_j$ is non-negative definite and $\alpha \subset \{1,2,\ldots,d\}$, then the quadratic form $\sum_{i,j \in \alpha} c_{i,j}x_ix_j$ is also non-negative definite. The quadratic form from the definition of the statistic $\tilde{T}_{SO}$ is non-negative definite (since it coincides with the non-negative definite quadratic form corresponding to the statistic $T_{SO}$). Therefore, the quadratic form on the right-hand side of \eqref{predstavlTSo} is also non-negative definite, and it follows from \eqref{equivTSoAndTildeTSo} and \eqref{predstavlTSo} that $T_{SO}$ coincides, up to $o_P(1)$, with the non-negative definite quadratic form (with coefficients independent of $n$) of the quantities $\mathcal{Z}'_{i_1,\ldots, i_m}$ for which the condition $\sigma_{i_1,\ldots,i_m} \ne 0$ holds. It was proved above that if $\sigma_{i_1,\ldots,i_m} \ne 0$, then the quantity $\frac{\mathcal{Z}'_{i_1,\ldots, i_m}}{\sigma_{i_1,\ldots,i_m}}$ is a summing statistic. Therefore, $T_{SO}$ coincides, up to $o_P(1)$, with the corresponding quadratic statistic.

\begin{remark}
We call a generalized quadratic statistic constructed from the summing statistics $T_{sum}^{[1]}, \ldots, T_{sum}^{[Q_{quad}]}$ a statistic $T_{gquad}$ that is a quadratic form (not necessarily non-negative definite) of the statistics $T_{sum}^{[1]}, \ldots, T_{sum}^{[Q_{quad}]}$. Quadratic statistics are a special case of generalized quadratic statistics. For a set of statistics $T_{sum}^{[1]}, T_{lb}^{[1]}, T_{sb}^{[1]}, \ldots, T_{sum}^{[Q]}, T_{lb}^{[Q]}, T_{sb}^{[Q]}, T_{gquad}^{[1]},\ldots, T_{gquad}^{[J]}$, one can prove an analogue of Theorem \ref{th2}, repeating its proof almost verbatim. Using this analogue of Theorem \ref{th2}, applying Theorem 3 \cite{Provost}, one can obtain a criterion for the asymptotic independence of vectors whose components are statistics of the form $T_{sum}, T_{lb}, T_{sb},T_{gquad}$. The resulting criterion will be more cumbersome than the criterion from Theorem \ref{thMainAsNez}.
\end{remark}

The author expresses his profound gratitude to A.M. Zubkov for constant attention.


{}
\end{document}